\documentclass[letterpaper, 10 pt, conference]{ieeeconf}
\IEEEoverridecommandlockouts        % This command is only
\usepackage{bm}
\let\labelindent\relax
\usepackage{enumitem}
\usepackage{cite}
\usepackage{amsmath}
\usepackage{amssymb,amsfonts}
\usepackage[ruled,linesnumbered]{algorithm2e}
\usepackage{textcomp}
\usepackage{graphicx,color}
\usepackage{mathrsfs}
\usepackage{subfigure}
\usepackage{url}
\usepackage{color}
\usepackage{dsfont}
\usepackage{bbm}
\usepackage{booktabs}
\usepackage{array}
\usepackage[table]{xcolor} 
\usepackage{yfonts}
\usepackage{tikz}
\usepackage{tcolorbox}
\usepackage{mathtools}
\usepackage{multicol}
\usepackage[colorlinks = True, linkcolor = blue, citecolor = blue]{hyperref}
\usepackage{accents}

\newtheorem{theorem}{Theorem}[section]
\newtheorem{lemma}[theorem]{Lemma}

\renewcommand*{\thesubsubsection}{\alph{subsubsection}}

\newcommand{\remend}{\relax\ifmmode\else\unskip\hfill\fi\hbox{$\bullet$}}

\newcommand{\real}{\mathbb{R}}

\newcommand{\E}{\operatorname{\mathbb{E}}}

\newcommand{\M}{\mathcal{M}}
\newcommand{\G}{\mathcal{G}}
\newcommand{\Hc}{\mathcal{H}}
\newcommand{\F}{\mathcal{F}}

\DeclareSymbolFont{bbold}{U}{bbold}{m}{n}
\DeclareSymbolFontAlphabet{\mathbbold}{bbold}

\newcommand{\R}{\mathbb{R}}

\newcommand\oprocendsymbol{\hbox{$\square$}}
\newcommand\oprocend{\relax\ifmmode\else\unskip\hfill\fi\oprocendsymbol}

\renewcommand{\theenumi}{(\roman{enumi})}

\newcommand*{\QEDA}{\hfill\ensuremath{\blacksquare}}%

\graphicspath{{img/}}

\begin{document}
\title{ \LARGE \bf Direct vs Indirect Methods for Behavior-based
  Attack Detection} \author{Darshan Gadginmath, Vishaal Krishnan, and
  Fabio Pasqualetti \thanks{This material is based upon work supported
    in part by awards ARO W911NF-20-2-0267, AFOSR-FA9550-20-1-0140 and
    AFOSR-FA9550-19-1-0235. Darshan Gadginmath and Fabio Pasqualetti
    are with the Department of Mechanical Engineering, University of
    California at Riverside, Riverside, CA, 92521, USA. Vishaal
    Krishnan is with the School of Engineering and Applied Sciences,
    Harvard University, Cambridge, MA, 02138, USA.  E-mail:
    \href{mailto:dgadg001@ucr.edu}{\texttt{dgadg001@ucr.edu}}, \href
    {mailto:vishaalk@ucr.edu}{\texttt{vkrishnan@seas.harvard.edu}},
    \href
    {mailto:fabiopas@engr.ucr.edu}{\texttt{fabiopas@engr.ucr.edu.}} }
}
\maketitle

\thispagestyle{empty}

\begin{abstract} We study the problem of data-driven attack detection for
 unknown LTI systems using only input-output behavioral data. In contrast with
 model-based detectors that use errors from an output predictor to detect
 attacks, we study behavior-based data-driven detectors. We construct a
 behavior-based chi-squared detector that uses a sequence of inputs and outputs
 and their covariance. The covariance of the behaviors is estimated using data
 by two methods. The first (direct) method employs the sample covariance as an
 estimate of the covariance of behaviors. The second (indirect) method uses a
 lower dimensional generative model identified from data to estimate the
 covariance of behaviors. We prove the consistency of the two methods of
 estimation and provide finite sample error bounds. Finally, we numerically
 compare the performance and establish a tradeoff between the methods at
 different regimes of the size of the data set and the length of the detection
 horizon. Our numerical study indicates that neither method is invariable
 superior, and reveals the existence of two regimes for the performance of the
 two methods, wherein the direct method is
 superior in cases with large data sets relative to the length of the detection horizon, while
 the indirect method is superior in cases with small data sets. 
\end{abstract}

% \begin{keywords}
% data-driven attack detection, direct and indirect estimation, system identification
% \end{keywords}

\section{Introduction} \label{sec:intro} 

Cyber-physical systems are growing in complexity and size since the advent of
better communication and higher computation power. This has also introduced
a greater possibility for adversarial attacks. Several attacks, such as the
ones on the power grid in Ukraine in 2015, the Maroochy attack in 2000, and
others mentioned in \cite{SMD-MP-etal:2019}, have exposed the vulnerability of
CPSs. Detection and mitigation of attacks has broadly been addressed using
model-based \cite
{FP-FD-FB:10y,CZB-FP-VG:2017,YM-RC-BS:2013,CM-JR:2016,RT-CM-JR:2018,YL-LS-TC:2017,AYL-GHY:2022}
and data-driven techniques \cite
{JG-SA-MT-ZSL:17,MK-AS:2021,YL-etal:2021,AT-etal:2014,VK-FP:2020,EN-KK:2017,MT-KK-IS-NM:2021}.
Model-based attack detection and their fundamental limitations have been
understood well \cite{FP-FD-FB:10y,CZB-FP-VG:2017}, model-based techniques assume
knowledge of the underlying system or the knowledge of the statistics of the
measured signals. Therefore the implementation of model-based methods need
prior system identification, which could be difficult to achieve due to
the complexity of the system or lack of data. In contrast, data-driven methods
usually operate with no knowledge of the system, yet offer scalable and
accurate detection. However, their limitations are not fully understood. We
study the problem of data-driven attack detection in stochastic systems and, in
particular, the tradeoffs between direct and indirect methods of detection. 

\textbf{Related work:} Machine Learning approaches to attack detection \cite
 {JG-SA-MT-ZSL:17,MK-AS:2021,YL-etal:2021} have gained popularity recently due
 to their ease of implementation. However, these techniques lack an insight
 into their functioning. System theoretic approaches provide more explainable
 solutions, and they can broadly be classified into direct and indirect
 methods. Direct methods map observations to nominal or attacked
 conditions \cite{AT-etal:2014} using past data. Indirect methods \cite
 {VK-FP:2020,EN-KK:2017,MT-KK-IS-NM:2021} utilize information about the
 underlying system to detect attacks. \cite{VK-FP:2020} considers the noiseless
 case where future outputs are predicted from past data and predictions are
 compared with measurements, while \cite{EN-KK:2017,MT-KK-IS-NM:2021} address
 the noisy case. \cite{EN-KK:2017} proposes a data-driven filter for fault
 detection and isolation. \cite{MT-KK-IS-NM:2021} analyses the conditions for
 successful attacks when the adversary does not have the knowledge of the
 system parameters. Differently from these works, this paper focuses on a finite-sample analysis
 and tradeoffs between direct and indirect methods.

 A promising model-based approach that can be adapted for the data-driven
 problem is the $\chi^2$ anomaly detector \cite
 {YM-RC-BS:2013,CM-JR:2016,RT-CM-JR:2018}. These works use the detector in
 conjunction with innovations from a state estimator and assume the knowledge
 of the system. We provide a modified data-driven approach using a $\chi^2$
 detector constructed from the measured system behaviors.

\textbf{Contributions:} This paper contributes a data-driven attack detector
 that uses the behaviors of an unknown LTI system to differentiate between
 nominal operation and attacked operation. To implement the detector, we propose two methods to
 estimate the covariance of the behaviors, namely the direct and indirect method. The
 direct method is motivated by discriminative classification where we use data
 to differentiate between nominal and attacked conditions. The
 indirect method is motivated by generative classification and the fact that an
 underlying lower dimensional system generates the observed behaviors. We analytically show that both methods are consistent. We
 also provide finite sample error bounds for the estimates of the two methods.
 Finally, we numerically compare the two methods at different noise levels
 using the ROC curve as a metric for performance, and establish different
 regimes where one outperforms the other.

\section{Data-driven attack detection} We consider the following
 discrete-time stochastic system:
\begin{align}
\label{eqn:system_dynamics}
\begin{aligned} 
    x_{t+1} &= Ax_t + Bu_t + w_t + B^a u^a_t, \\    
    y_t &= Cx_t + v_t + G^a y^a_t,
\end{aligned}   
\end{align} where $A \in \mathbb{R}^{n\times n}, B \in \mathbb{R}^{n\times m},
 C \in \mathbb{R}^{p\times n} $ are the system matrices, $x_t \in \real^n$,
 $u_t \in \real^m$, and $y_t \in \real^p$ are the state, control input and
 the output, respectively, $u^a_t$ and $y^a_t$ are the malicious input and false
 sensor measurement injected by an adversary as attacks on the system. The attacks are
 affecting the system through the matrices $B^a$ and $G^a$, respectively. We assume that the system is stable, controllable, and observable. The process noise $w_t$ and measurement noise $v_t$ are i.i.d.
 Gaussian processes with zero mean. The data-driven attack detection
 problem is posed as follows.

{\bf Data-driven attack detection:} Given inputs $\bm{u} = \left
[u_0^\top, \dots, u_{T-1}^\top \right]^\top$ and outputs $\bm{y} = \left
[y_1^\top, \dots, y_T^\top \right]^\top$ over the detection horizon $T$,
determine if $(\bm{u^a},\bm{y^a}) \not\equiv 0 $, where $\bm{u^a}  =[
{u^a_0}^\top, \dots, {u^a_{T-1}}^\top ]^\top$ and $\bm{y^a} = [
{y^a_1}^\top, \dots, {y^a_T}^\top]^\top$. 

To perform attack detection, we have input-output data from $N$ attack-free
experiments $\left(\bm{u}^{(i)},\bm{y}^{(i)}\right)$, where $i \in \{1,2,\dots,N\}$. Each experiment runs over a horizon of length $T$, with $x_0 = 0$. The inputs for the experiments are generated
by a Gaussian random process $\mathcal{N}(0,\Sigma_u)$. 

 We begin by noting that the outputs $\bm{y}$ are generated by inputs $\bm{u}$, noises $\bm{w}$ and $\bm
 {v}$, and attacks $\bm{u}^a$ and $\bm{y}^a$ as,
 \begin{align}
 \bm{y} = \mathcal{C} \bm{u} + \mathcal{C}' \bm{w} + \bm{v} + \mathcal{C}^a \bm{u^a} + F^a \bm{y^a}. 
 \end{align}
Here $\bm{w}$ and $\bm{v}$ are the process and measurement noise $w_t$ and $v_t$
over the horizon $T$. The matrix $\mathcal{C}$ is defined as 
\begin{align}
\mathcal{C} = \begin{bmatrix}
CB & 0 & \dots & 0 \\
CAB & CB & \dots & 0 \\ 
\vdots & \vdots & \ddots & \vdots \\
CA^{T-1}B & CA^{T-2}B & \dots & CB 
\end{bmatrix}.
\end{align}
$\mathcal{C}'$ is obtained by replacing $B$ with the identity matrix $I$ in $\mathcal{C}$, and
$\mathcal{C}^a$ is obtained by replacing $B$ with $B^a$ in $\mathcal{C}$. Note
that in the attack-free case, the assumptions made earlier ensure that $\bm
{y}$ as well as the behavior $Z$ are stationary Gaussian processes. In general, control inputs need not be generated by a Gaussian distribution. However, any output stabilizing controller will generate stationary measurements when the system is in operation.

Although we have access to the inputs $\bm{u}$ and outputs $\bm{y}$ over the
detection horizon, we have access to neither the system matrices nor the Markov
parameters $CA^iB$, making it difficult to perform attack detection using
methods that rely on a model based $\chi^2$ detector or KL-divergence \cite{CZB-FP-VG:2017}.
We propose to modify the usage of the $\chi^2$ detector as follows,
\begin{align} g = Z^\top S^{-1} Z \overset{H_0}{\underset{H_1}{\lessgtr}} \lambda,
\end{align} 
where $Z =[\bm{u}^\top,\bm{y}^\top]^\top$ is the behavior of the
 system \eqref{eqn:system_dynamics} over the time horizon $T$. The
 covariance matrix $S = \E[ZZ^\top]$. Note that $\E[Z] = 0$ under nominal
 operation. Further, hypothesis $H_1$ denotes an alarm, whereas $H_0$ denotes
 nominal operation. $\lambda$ is an arbitrary threshold used to differentiate
 nominal and attacked operation. To implement the above detector, we need to
 estimate $S$ for which we propose the direct and indirect methods.

\subsection{Direct method} The direct method is motivated by discriminative
 methods of classification, where a measured signal is mapped to nominal or
 attacked operation using data. We propose to use the data from $N$ experiments to
 estimate the covariance $S$ without considering the structure or
 characteristics of the samples. In particular, we use the sample covariance as
 a direct estimate,
\begin{align}
    \hat{S}^d &= \frac{1}{N}\sum_{i=0}^N \bm{z}^{(i)} {\bm{z}^{(i)}}^\top = \frac{1}{N} \bm{Z}\bm{Z}^{\top}. 
\end{align} 
Here, $\bm{z}^{(i)} = [{u^{(i)}_0}^\top,\dots,{u^{(i)}_{T-1}}^\top,{y^{
(i)}_1}^\top,\dots,{y^{(i)}_T}^\top]^\top \in \mathbb{R}^{T(m+p) \times 1}$ ,
and $\bm{Z} = [\bm{z}^{(1)} \ \bm{z}^{(2)} \ \dots \ \bm{z}^{(N)}]  \in \mathbb
{R}^{T(m+p) \times N}$. The direct method is computationally simple. In cases
with large data sets, the sample covariance is a good estimator of the true
covariance matrix $S$. As the dimension of the covariance matrix $S$ grows with the detection horizon $T$, the
sample covariance needs more data for accurate estimation.
Further, when $N < T(m+p)$, the sample covariance is not a good estimator.
Therefore, in the sequel, we propose the indirect method which can potentially
outperform the direct method in cases with little data. 

\subsection{Indirect method} 

In the indirect method, we identify a lower dimensional generative model that
gives rise to the behaviors. We first break the behavior $Z$ into smaller
behaviors called \emph{minor behaviors}, each of length $L$ denoted by $f_t$.
In particular,  $f_t = [u_{t-L}^\top \ u_{t-L+1}^\top \ \dots \ u_{t}^\top | y_
{t-L}^\top \ y_{t-L+1}^\top \ \dots y_{t}^\top]^\top$. Note that we can
construct $T-L+1$ minor behaviors from $Z$. We seek to regress future minor
behaviors on past minor behaviors of the form, $f_{t+1} = \M f_t + \epsilon_t$
using data. $\M$ defines the linear evolution of the minor behaviors while
$\epsilon_t$ is additive noise entering the system which is uncorrelated to
the behaviors. By exploiting the nature of the underlying dynamics of the minor
behaviors, we can potentially estimate the covariance matrix of the behaviors
with less data. We first note that the minor behaviors follow a stationary
process owing to our initial assumptions. Let the covariance of the minor
behaviors be $P$, and the covariance of the noise $\epsilon_i$ be
$\Sigma_\epsilon$. Now, the covariance of the behaviors can be computed as,
\begin{align}
\label{eqn:indirect-lyap}
\begin{aligned}
\mathrm{Cov}\left(f_{i+1}\right) &= \mathrm{Cov}\left(\M f_i + \epsilon_i\right), \\
P &= \M P {\M}^\top + \Sigma_\epsilon. 
\end{aligned}
\end{align}  
It is clear that equation \eqref{eqn:indirect-lyap} is in the form of a
Lyapunov equation. Therefore, the covariance matrix $P$ can be computed
from \eqref{eqn:indirect-lyap} if $\M$ and $\Sigma_\epsilon$ are available. 
If $\M$ is regressed from data, $\Sigma_\epsilon$ can be estimated from
the residuals of the regression using the sample covariance. In turn, $P$ can
be estimated by solving the Lyapunov equations using the estimates of $\M$ and
$\Sigma_\epsilon$. We propose to estimate $\M$ using Ordinary Least Squares
(OLS) as,
\begin{align}
\label{eqn:indirect-regression}
\hat{\M} &= F'F^\top\left(FF^\top\right)^{-1},
\end{align} 
where $F$ and $F'$ are data matrices constructed from the minor behaviors
obtained from the experimental data:
\begin{align*}
F &= \left[f^{(1)}_1 \ \dots \ f^{(1)}_{T-L} \ f^{(2)}_{1}\ \dots \ f^{(N)}_{T-L}\right] \in \R^{L(p+m) \times N(T-L)}  , \\
F' &= \left[f^{(1)}_2 \ \dots \ f^{(1)}_{T-L+1} \ f^{(2)}_{2}\ \dots \ f^{(N)}_{T-L+1}\right]. \label{eqn:indirect-data}
\end{align*}
If $\hat{\epsilon}_i$ are the residuals from the regression \eqref{eqn:indirect-regression},
$\Sigma_\epsilon$ can be estimated using the sample covariance. 
\begin{align}
\hat{\epsilon}^{(i)}_j &= f^{(i)}_{j+1} - \hat{\M} f^{(i)}_j,  \\
\hat{\Sigma}_\epsilon &= \frac{1}{N_{id}}EE^\top.
\end{align}
Here $E = [\hat{\epsilon}^{(1)}_1 \ \hat{\epsilon}^{(1)}_2 \dots \ \hat{\epsilon}^{(N)}_{T-L}]$ and $N_{id} = N(T-L)$. $f^{(i)}_j$ represents the minor behavior
$f_j$ from the experiment $i$ with $j \in \{1,2,\dots,T-L+1\}$. We now have the estimates of $\M$ and $\Sigma_\epsilon$ to solve equation \eqref{eqn:indirect-lyap}. Therefore, $P$ is estimated by solving $\hat{P}=\hat{M}\hat{P}\hat{M}+\hat{\Sigma}_\epsilon$.

We next describe the estimation of $S$ using the covariance of the minor behaviors $P$. Let $D
= \left[f_1^\top \ f_2^\top \ \dots \ f_{T-L+1}^\top \right]^\top$. Then $Z =
KD$, where $K$ is a known sparse matrix that reconstructs $Z$ from $D$. Let the covariance of $D$
be $\E\left[DD^\top \right] = \Sigma_D$. Therefore $S = \E[ZZ^\top] =  K\Sigma_{D}K^\top$. Let the matrix $\F$ be given as,
\begin{align}
\F&= 
\begin{bmatrix}
I & \M & \dots & \M^{T-L} \\
\M & I & \dots & \M^{T-L-1} \\
\vdots &\vdots &\ddots &\vdots \\
\M^{T-L} & \M^{T-L-1} & \dots & I
\end{bmatrix}. \nonumber
\end{align}
Then the covariance $\Sigma_D$ is,
\begin{align}
\Sigma_D &= \F \otimes P,
\label{eqn:indirect-sig-est}
\end{align}
where $\otimes$ denotes the Kronecker product.
Therefore, $\Sigma_D$ can be estimated as $\hat{\Sigma}_D = \hat
{\F} \otimes \hat{P}$, where $\hat{\F}$ is obtained by replacing $\M$ with
its estimate $\hat{\M}$. Finally, $S$ is estimated as $\hat{S}^{id} = K\hat{\Sigma}_DK^\top$.

% The indirect method utilizes underlying information of the system
% which could potentially make it a better estimator of the covariance matrix in
% low data regimes. In the following section, we establish the consistency as
% well finite-sample analysis on the deviation of the estimate of the indirect
% method. 
% \begin{algorithm}
% \label{algo:indirect}
% \caption{Indirect method to estimate $S$
% }
% {\bf{\emph{Data}}} $:= (\bm{u}^{(i)},\bm{y}^{(i)})$, $\forall i \in \{1,\dots, N\}$ \\
% {\bf{\emph{Output}}} $:= \hat{S}^{id}$\\
% \label{algo:indirect}
% \SetAlgoLined
% \txt{Construct data matrices $F$ and $F'$ \eqref{eqn:indirect-data}} \\
% \txt{Solve OLS problem, } $\hat{\M} = F'F^\top \left(FF^\top\right)^{-1}$ \\
% \txt{Compute residuals from OLS solution} $\epsilon_i = f_{i+1} - \hat{\M}f_i$ , $\forall i \in \{1,\dots,N_{id}\}$ \\
% \txt{Construct residual matrix}  $E = \left[\hat{\epsilon}_1 \ \hat{\epsilon}_2 \ \dots \ \hat{\epsilon}_{N_{id}} \right]$. \\
% \txt{Estimate $\Sigma_\epsilon$ as $\hat{\Sigma}_\epsilon = \frac{1}{N_{id}}EE^\top $}. \\
% \txt{Solve $\hat{P} = \hat{\M} \hat{P} \hat{\M}^\top + \hat{\Sigma}_\epsilon$ to obtain $\hat{P}$}. \\
% \txt{Compute $\hat{\Sigma}_D$ as in \eqref{eqn:indirect-sig-est} using $\hat{\M}$ and $\hat{P}$}. \\
% \txt{Obtain estimate of $S$ as $\hat{S} = K\hat{\Sigma}_D K^\top$}
% \end{algorithm}

\section{Consistency and Finite-Sample Analysis}

In this section, we establish that both the direct and indirect methods are consistent
estimators of the covariance matrix $S$. We also provide finite sample analysis
on the deviation of the estimated covariance matrices $\hat{S}^d$ and $\hat{S}^
{id}$ from the true covariance.

\begin{theorem}\label{thm:consistency}
{\bf (\emph{Consistency})}
Let the covariance matrix $S$ be estimated as $\hat{S}^d$ and $\hat{S}^{id}$ using the direct and indirect methods, respectively, then,
\begin{align*}
\lim_{N \rightarrow \infty} \hat{S}^d = \lim_{N \rightarrow \infty} \hat{S}^{id}=S.
\end{align*}
\oprocend
\end{theorem}
We refer the reader to Appendix A for the proof. This result establishes that with infinite data, both the methods estimate the covariance matrix $S$ perfectly.
We next establish a finite sample error bound for the two methods. We denote the spectral norm of a matrix by $\|\cdot\|$.
\begin{theorem}
{\bf (\emph{Finite-sample bound on error - Direct method})}
\label{thm:finitesample-direct} 
If the covariance $S$ is estimated as $\hat{S}^d$ using the direct method, then
for any $\theta \geq 0$, and $r = \frac{\mathrm{Tr}(S)}{\|S\|}$,
\begin{align}
\|S - \hat{S}^d\| \leq \left(\sqrt{\frac{2\theta (r+1)}{N}} + \frac{2\theta r}{N}\right) \| S \| \label{eqn:thm-direct}
\end{align}
with probability at least $1 - 2  T(m+p) e^{-\theta}$, for any $ N > T(m+p)$. \hfill $\square$
\end{theorem}
We provide the proof in Appendix B. Theorem \ref{thm:finitesample-direct} characterizes a finite sample
 bound on the error in estimation arising from
 using the direct method. This bound grows larger with the dimension of $S$ as
 $r$ is a strictly increasing function of the size of $S$. This contributes to
 high sample complexity of the direct method, leading to poor estimates at low
 data regimes for big $S$ matrices. Further, the above bound converges to 0 at
 a rate of the order $\mathcal{O}(\frac{1}{\sqrt{N}})$. 

For the indirect method, it is evident from equation \eqref
{eqn:indirect-sig-est} that solving the OLS problem introduces errors in
the matrices $\F$ and $P$. Define the errors arising from the OLS solution as
$\Delta_\M = \M - \hat{\M}$ and $\Delta_
{\Sigma_\epsilon} = \Sigma_\epsilon - \hat{\Sigma}_\epsilon$. Using
result from Theorem 1 of \cite{TS-AR:2019}, we obtain with probability at least $1-\theta$, 
\begin{align}
\label{eqn:ols-bnd}
\|\Delta_\M\| &\leq \sqrt{\frac{k}{N_{id}}} \ \gamma_s\left(\M,\frac{\theta}{4}\right),
\end{align}
for any $N_{id} \geq \max\Big(N_\eta(\theta), N_s(\theta) \Big)$. 
Here $k$ is an absolute constant. We define $\gamma_s(\M,\frac{\theta}
{4})$, $N_\eta(\theta)$ and $N_s(\theta)$ as follows,
\begin{align*}
\gamma_s(\M,\theta) &= \sqrt{8L(m+p) \Big(\log\frac{5}{\theta}+ \frac{\log 4 \mathrm{Tr}(\Gamma_N(\M))+1}{2} \Big) }, \\
 N_\eta(\theta) &= k \log\frac{2}{\theta}+(L(m+p)\log 5), \\
 N_s(\theta) &= k \Big(L(m+p)\log(\mathrm{Tr}(\Gamma_N(\M))+1) \\ & \quad+2L(m+p)\log\frac{5}{\theta}\Big), \\
\Gamma_N(\M) &= \sum_{j=0}^{N_{id}} \M^j {\M^j}^\top.
\end{align*}
Equation \eqref{eqn:ols-bnd} is a tight probabilistic bound on the error arising
from the OLS problem. The error in covariance of the residuals $\Delta_
{\Sigma_\epsilon}$ follows the same distribution as established in
equation \eqref{eqn:thm-direct} because the residuals of the OLS solution are
uncorrelated Gaussian random variables. Therefore both $\|\Delta_\M\|$ and $\|\Delta_
{\Sigma_\epsilon}\|$ converge to $0$ with rate $\mathcal{O}\left(\frac{1}{\sqrt
{N}}\right)$.

Before we bound the error in estimation from the indirect
method, we first provide a sensitivity analysis for matrices $\F$ and $P$ when
they are estimated. 
\begin{lemma}{\bf \emph{(Sensitivity of $\F$ and $P$)}}
\label{lem:fp-bnd}
Given estimated matrices $\hat{\F}$ and $\hat{P}$, let the errors in estimation of $P$ be $\Delta_P = \hat{P}-P$, and $\F$ be $\Delta_\F = \hat{\F}-\F$. Then,
\begin{enumerate}[label = (\alph*)]
    \item Sensitivity of $P$ is,
\begin{align}
\|\Delta_P\| &\leq \left(\sqrt{L(m+p)}\|I \otimes I - \M^\top \otimes \M^\top\|\right) \nonumber \\ 
& \quad \Bigg((1+\|\M+\Delta_\M\|)^2\Big(\frac{\|\Delta_{\Sigma_\epsilon}\|}{\|\Sigma_\epsilon + \Delta_{\Sigma_\epsilon}\|}\Big) \nonumber \\ & \quad+ 2(\|\M\| + \|\Delta_\M\|)^2  \Big(\frac{\|\Delta_{\M}\|}{\|\M+\Delta_\M\|}\Big) \Bigg).
\end{align}
if $\Delta_P$ and $\Delta_{\Sigma_\epsilon}$ are positive semi-definite. 
\item Sensitivity of $\F$ is, 
\begin{align}
\|\Delta_{\F}\| &\leq \|\F\| + 1 + 2 \big(\|\M+\Delta_\M\| \nonumber \\
&\quad + \|\M+\Delta_\M\|^2 + \dots + \|\M+\Delta_\M\|^{T-L}\big).
\end{align}
\end{enumerate}
\oprocend
\end{lemma}
We provide the proof in Appendix C. The above sensitivity analysis allows us to provide the rate of convergence for both $\Delta_P$ and $\Delta_\F$. Since $\|\Delta_\M\|$ and $\|\Delta_{\Sigma_\epsilon}\|$ both converge at the rate of $\mathcal{O}\left(\frac{1}{\sqrt{N}}\right)$, $\|\Delta_P\|$ converges to 0 at the rate of $\mathcal{O}\left(\frac{1}{\sqrt{N}}\right)$, and $\|\Delta_\F\|$ converges to a system dependent constant at the rate of $\mathcal{O}\left(\frac{1}{\sqrt{N}}\right)$.

Following Lemma \ref{lem:fp-bnd}, we provide a probabilistic bound on the
estimation error of the indirect method obtained from finite samples. We bound
the estimation error $\|S-\hat{S}^{id}\|$ as a function of $\Delta_P$ and
$\Delta_{\F}$. Recall that $S = K\Sigma_DK^\top$ and $\hat{S}^{id} = K \hat
{\Sigma}_D K^\top$. Here, the matrix $K$ has at most one identity matrix per
column as it picks elements from the vector $D$ to construct $Z$, therefore
$\|K\| = 1$. Using the properties of the Kronecker product and the triangle
inequality, we obtain the error bound.
\begin{align}
\|S-\hat{S}^{id}\| & = \|K(F\otimes P)K^\top \nonumber\\ & \quad - K\left((F+\Delta_F) \otimes (P+\Delta_P)\right)K^\top\| \nonumber \\ 
&\leq \|F\otimes P - (F+\Delta_F) \otimes (P+\Delta_P)\| \nonumber \\
&\leq  \|\F\| \ \|\Delta_P\| + \|\Delta_\F\| \ \|P\| + \|\Delta_\F\| \ \|\Delta_P\|. \label{eqn:indirect-finite-sample}
\end{align} The indirect method scales at the rate $\mathcal{O}(L)$ while the
 direct method scales at the rate $\mathcal{O}(T)$. The data set expands for
 the indirect method to $N(T-L)$ as each experiment $\bm{z}^{(i)}$ is broken
 into minor behaviors $f^(i)_j$ with $j \in \{1,\dots,T-L+1\}$. This makes the
 indirect method better suited for cases with low data sizes and longer
 detection horizons as $L \leq T$. The performance of the indirect method is also dependent
 on the choice of $L$. When $L=n$, the representation of the minor behaviors in
 equation \eqref{eqn:indirect-lyap} is exact as proven in Appendix A. However,
 choosing a large value for $L$ increases the model complexity, thereby
 overfitting the data. Contrary to this, choosing a small $L$ can underfit the
 data. The direct method is free of design parameters and can potentially
 outperform the indirect method when the size of the available data set $N$ is
 large. 

\section{Simulations} 

In this section we numerically compare and establish a tradeoff in performance
of the detector using the two proposed methods of estimation. To make a
comparison, we consider the case when the system is operating nominally as well
as when the system is attacked by an adversary. As the performance of the
detector depends on both nominal and attacked conditions, we compare the False
Positive Rate (FPR) and True Positive Rate (TPR) of the two methods. We also
investigate the effect of the noise on the performance of the two
methods by varying the signal to noise ratio. 

We consider a stable SISO system with $n=3$ with randomly generated $A,B$, and
$C$ matrices. Also, for experimental data we use $u_t \sim \mathcal{N}
(0,\sigma_u I)$, $w_t \sim \mathcal{N}(0,\sigma_w I)$, and $v_t \sim \mathcal
{N}(0,\sigma_v I)$. We vary $\sigma_u$, $\sigma_v$ and $\sigma_w$ to compare
the performance of the estimates under various signal to noise ratios using ROC
curves. 

Under nominal conditions, the performance of the detector is determined by the
FPR, which is defined as $\mathbb{P}(g \geq \lambda | (\bm
{u}^a,\bm{y}^a) \equiv 0)$. Under an attack, the performance is determined by
the TPR defined as $\mathbb{P}
(g \geq \lambda | (\bm{u}^a,\bm{y}^a) \not\equiv 0)$. To compute the TPR, we introduce a
detectable attack where $u^a_t \sim \mathcal{N}(0,1.5)$, and
$y^a_t \sim \mathcal{N}(0,1.5)$. The TPR and FPR in the following comparisons
have been averaged over 50 trials. For the indirect method, $L=n$.

\subsubsection{Comparison 1} 
\begin{figure*}
\centering
\includegraphics[width=0.30\textwidth]{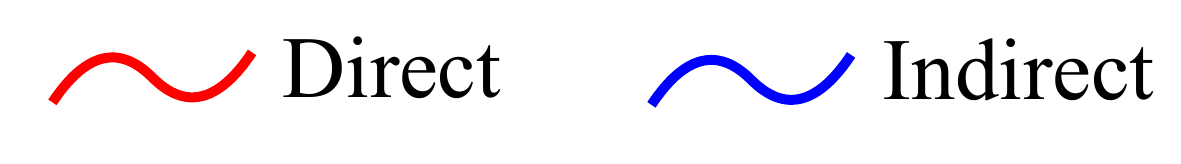}
\vspace*{-0.4cm}
\begin{multicols}{4}
\hspace*{0.75cm}
\begin{tikzpicture}
  \node (img1)  {\includegraphics[width=0.20\textwidth]{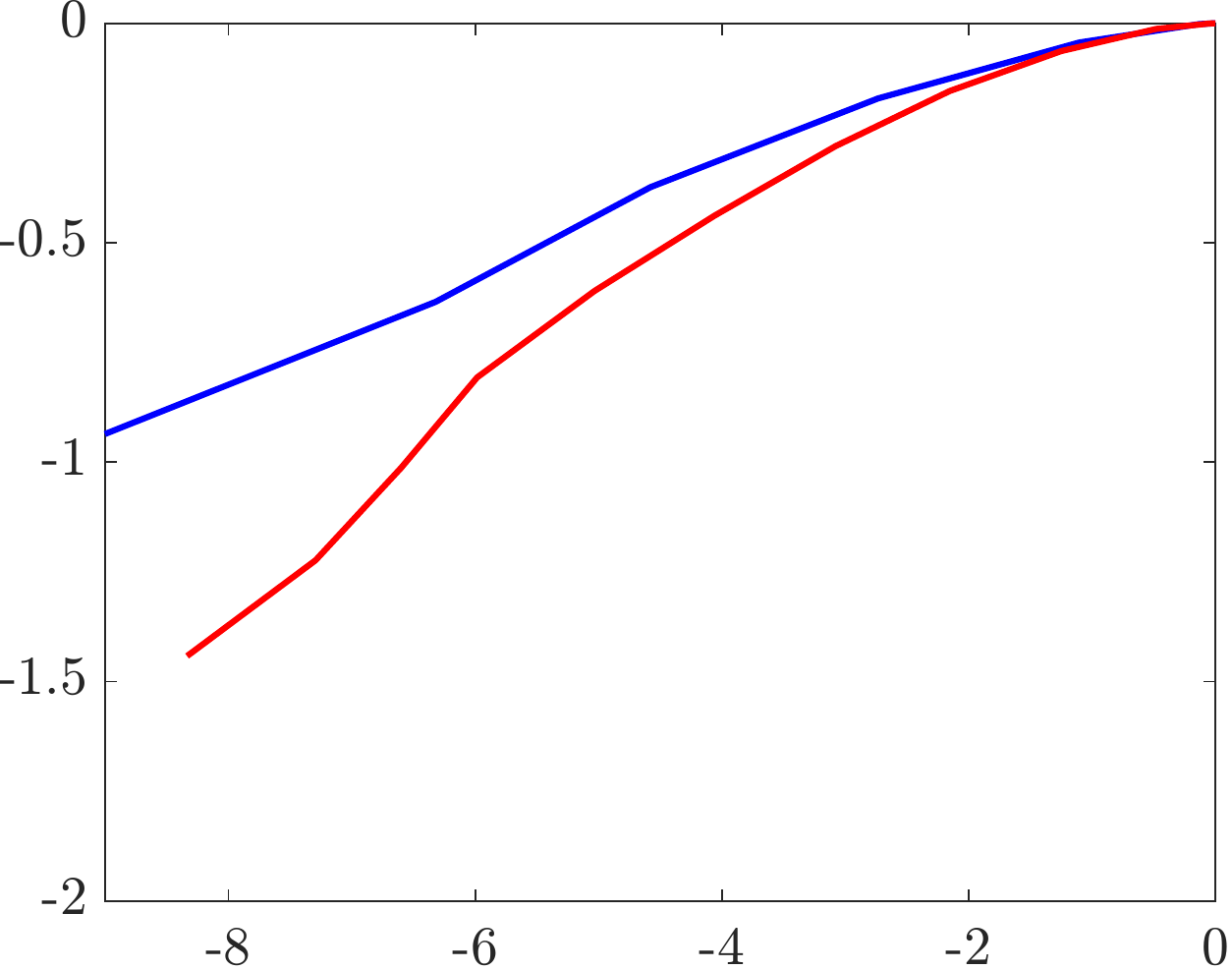}};
  \node[above of= img1, node distance=0cm, yshift=1.7cm,font=\color{black}]  {\small $N=40$};  
  \node[below of= img1, node distance=0cm, yshift=-1.7cm,font=\color{black}]  {\small $\log$ FPR};
  \node[above of= img1, node distance=0cm, xshift=-1.5cm,yshift=1.7cm,font=\color{black}]  {\small (a)};
  \node[left of= img1, node distance=0cm, rotate=90, anchor=center,yshift=2.1cm,font=\color{black}] { \small $\log$ TPR};
\end{tikzpicture}\columnbreak
\begin{tikzpicture}
  \node (img1)  {\includegraphics[width=0.20\textwidth]{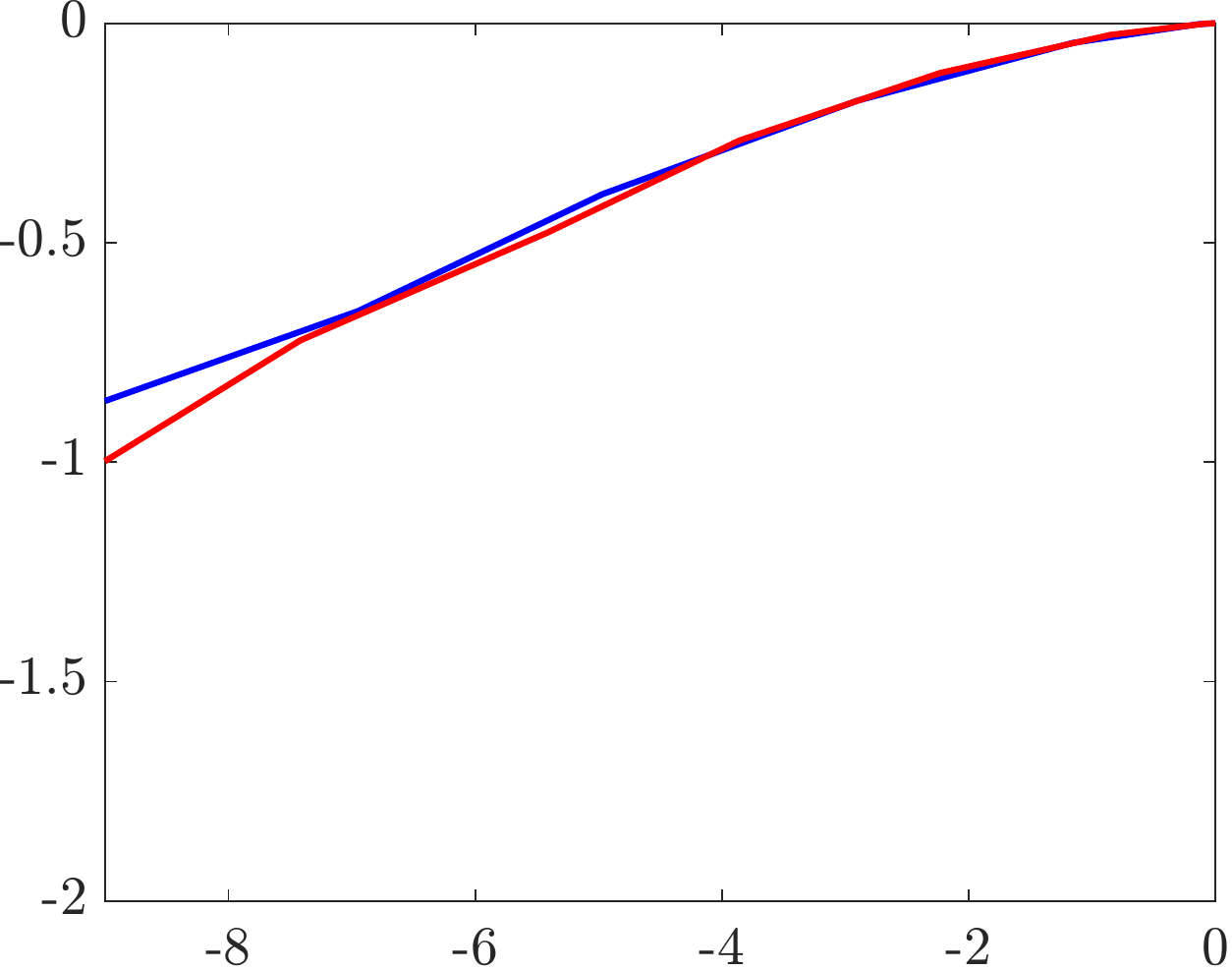}};
  \node[above of= img1, node distance=0cm, yshift=1.7cm,font=\color{black}]  {\small $N=90$}; 
  \node[below of= img1, node distance=0cm, yshift=-1.7cm,font=\color{black}]  {\small $\log$ FPR};
  \node[above of= img1, node distance=0cm, xshift=-1.5cm,yshift=1.7cm,font=\color{black}]  {\small (b)};
  % \node[left of= img1, node distance=0cm, rotate=90, anchor=center,yshift=2.1cm,font=\color{black}] { \small TPR};
\end{tikzpicture}\columnbreak 
\begin{tikzpicture}
  \node (img1)  {\includegraphics[width=0.20\textwidth]{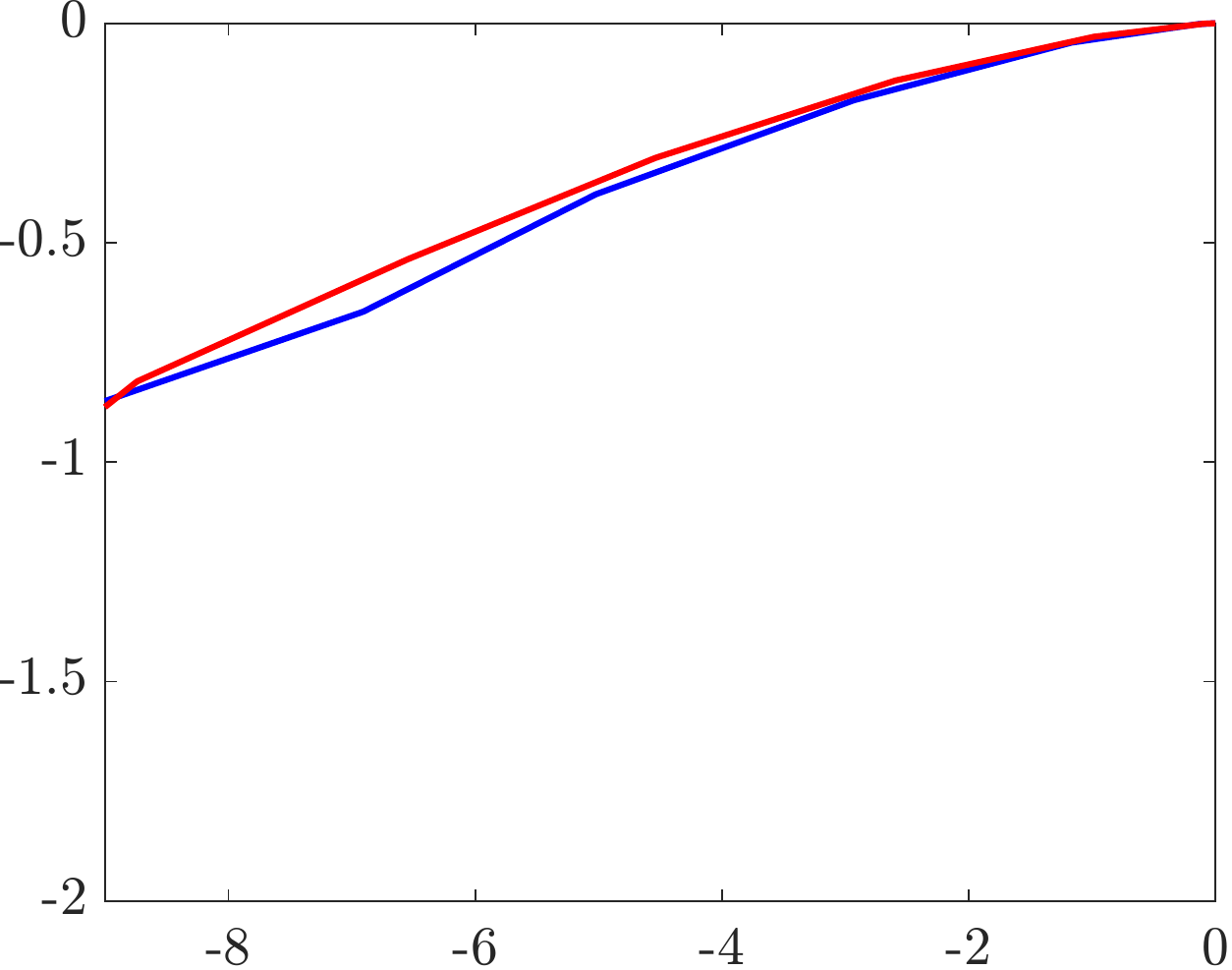}};
  \node[above of= img1, node distance=0cm, yshift=1.7cm,font=\color{black}]  {\small $N=150$};   
   \node[below of= img1, node distance=0cm, yshift=-1.7cm,font=\color{black}]  {\small $\log$  FPR};
  \node[above of= img1, node distance=0cm, xshift=-1.5cm,yshift=1.7cm,font=\color{black}]  {\small (c)};
  % \node[left of= img1, node distance=0cm, rotate=90, anchor=center,yshift=2.1cm,font=\color{black}] { \small TPR};
\end{tikzpicture}\columnbreak
\begin{tikzpicture}
  \node (img1)  {\includegraphics[width=0.20\textwidth]{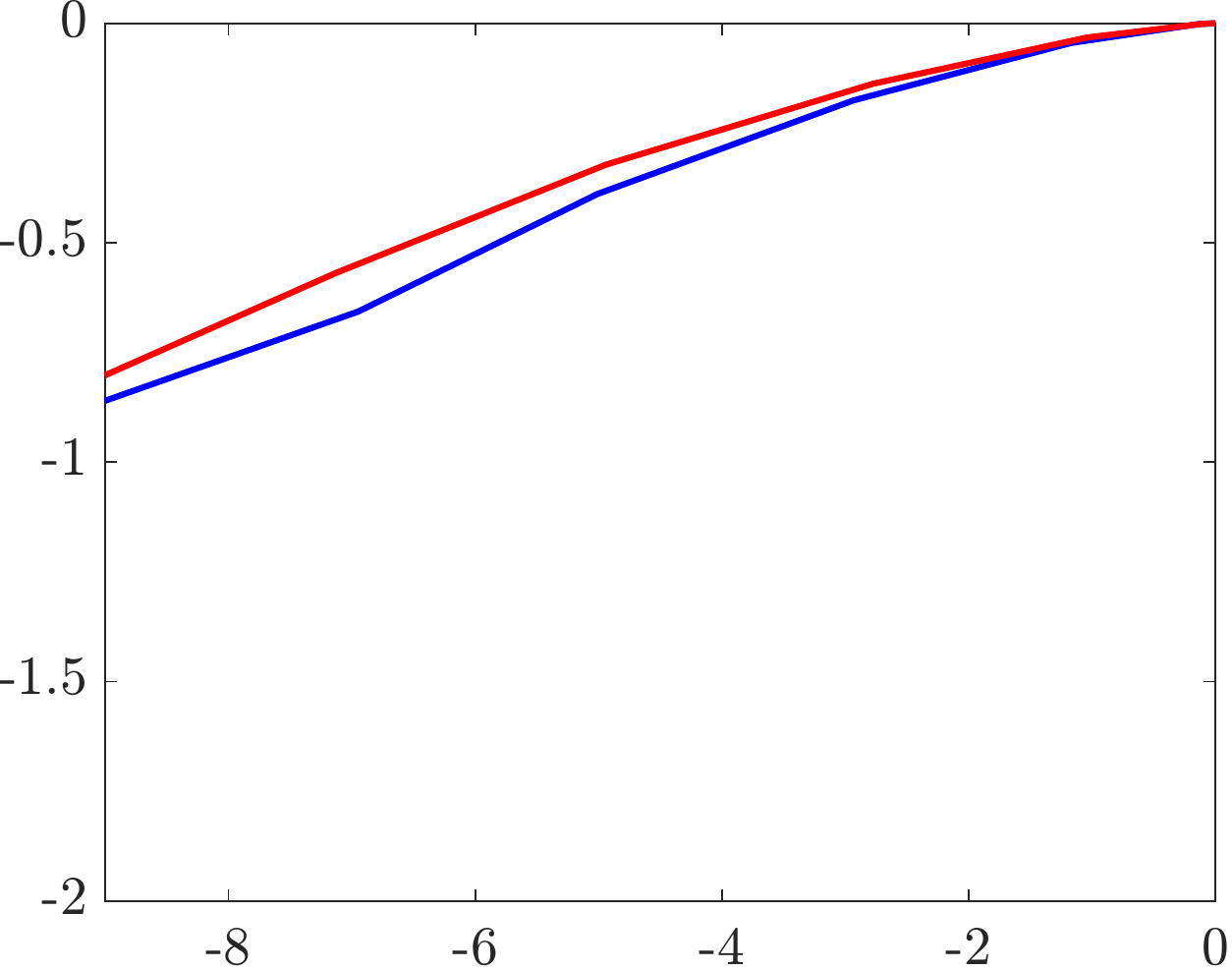}};
  \node[above of= img1, node distance=0cm, yshift=1.7cm,font=\color{black}]  {\small $N=200$};   
  \node[below of= img1, node distance=0cm, yshift=-1.7cm,font=\color{black}]  {\small $\log$  FPR};
  \node[above of= img1, node distance=0cm, xshift=-1.5cm,yshift=1.7cm,font=\color{black}]  {\small (d)};
  % \node[left of= img1, node distance=0cm, rotate=90, anchor=center,yshift=2.1cm,font=\color{black}] { \small TPR};
\end{tikzpicture}
\end{multicols}
\vspace*{-0.8cm}
\begin{multicols}{4}
\hspace*{0.75cm}
\begin{tikzpicture}
  \node (img1)  {\includegraphics[width=0.20\textwidth]{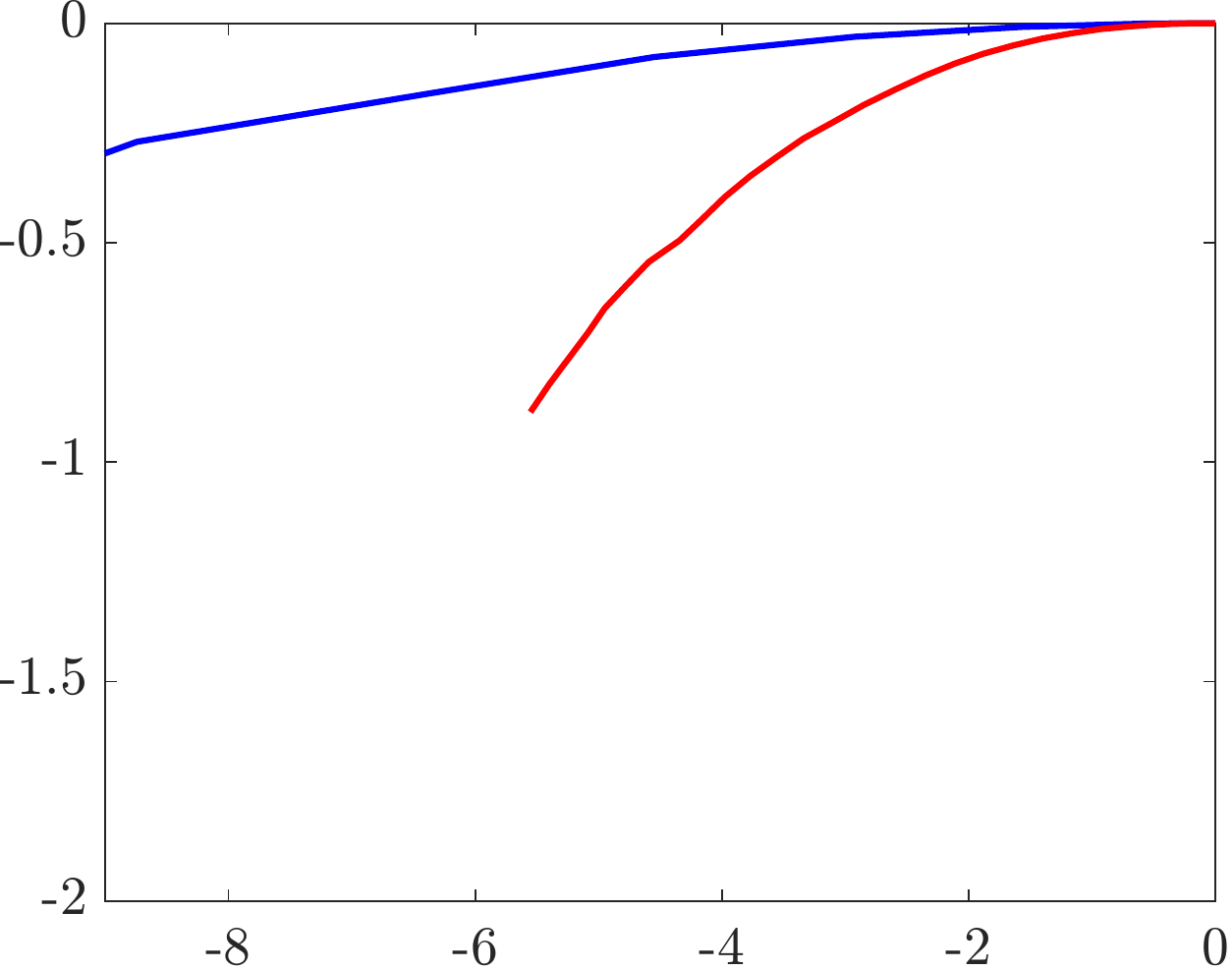}};
  \node[above of= img1, node distance=0cm, yshift=1.7cm,font=\color{black}]  {\small $N=40$};  
  \node[below of= img1, node distance=0cm, yshift=-1.7cm,font=\color{black}]  {\small $\log$ FPR};
  \node[above of= img1, node distance=0cm, xshift=-1.5cm,yshift=1.7cm,font=\color{black}]  {\small (e)};
  \node[left of= img1, node distance=0cm, rotate=90, anchor=center,yshift=2.1cm,font=\color{black}] { \small $\log$ TPR};
\end{tikzpicture}\columnbreak
\begin{tikzpicture}
  \node (img1)  {\includegraphics[width=0.20\textwidth]{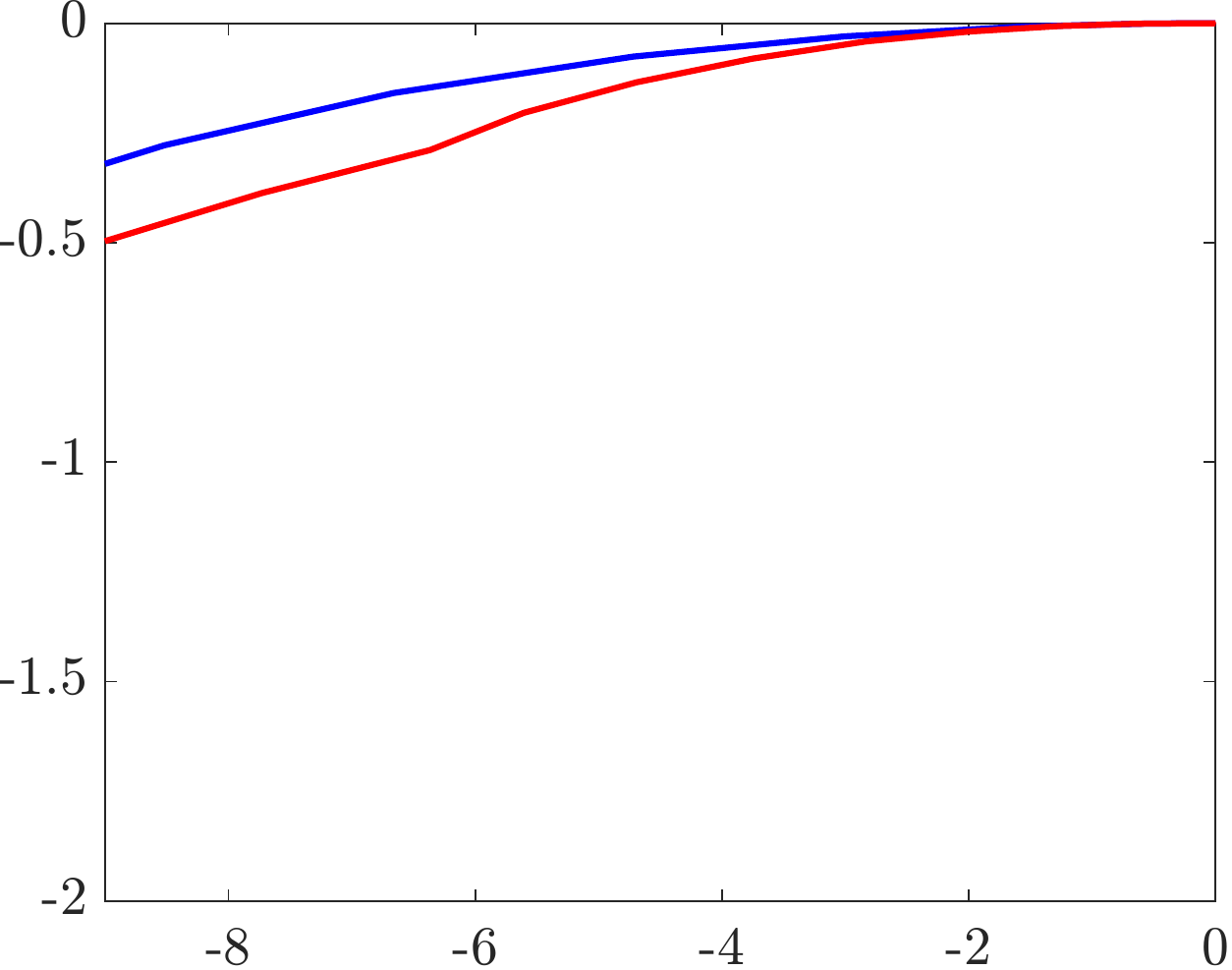}};
  \node[above of= img1, node distance=0cm, yshift=1.7cm,font=\color{black}]  {\small $N=90$}; 
   \node[below of= img1, node distance=0cm, yshift=-1.7cm,font=\color{black}]  {\small  $\log$ FPR};
    \node[above of= img1, node distance=0cm, xshift=-1.5cm,yshift=1.7cm,font=\color{black}]  {\small (f)};
  % \node[left of= img1, node distance=0cm, rotate=90, anchor=center,yshift=2.1cm,font=\color{black}] { \small TPR};
\end{tikzpicture}\columnbreak 
\begin{tikzpicture}
  \node (img1)  {\includegraphics[width=0.20\textwidth]{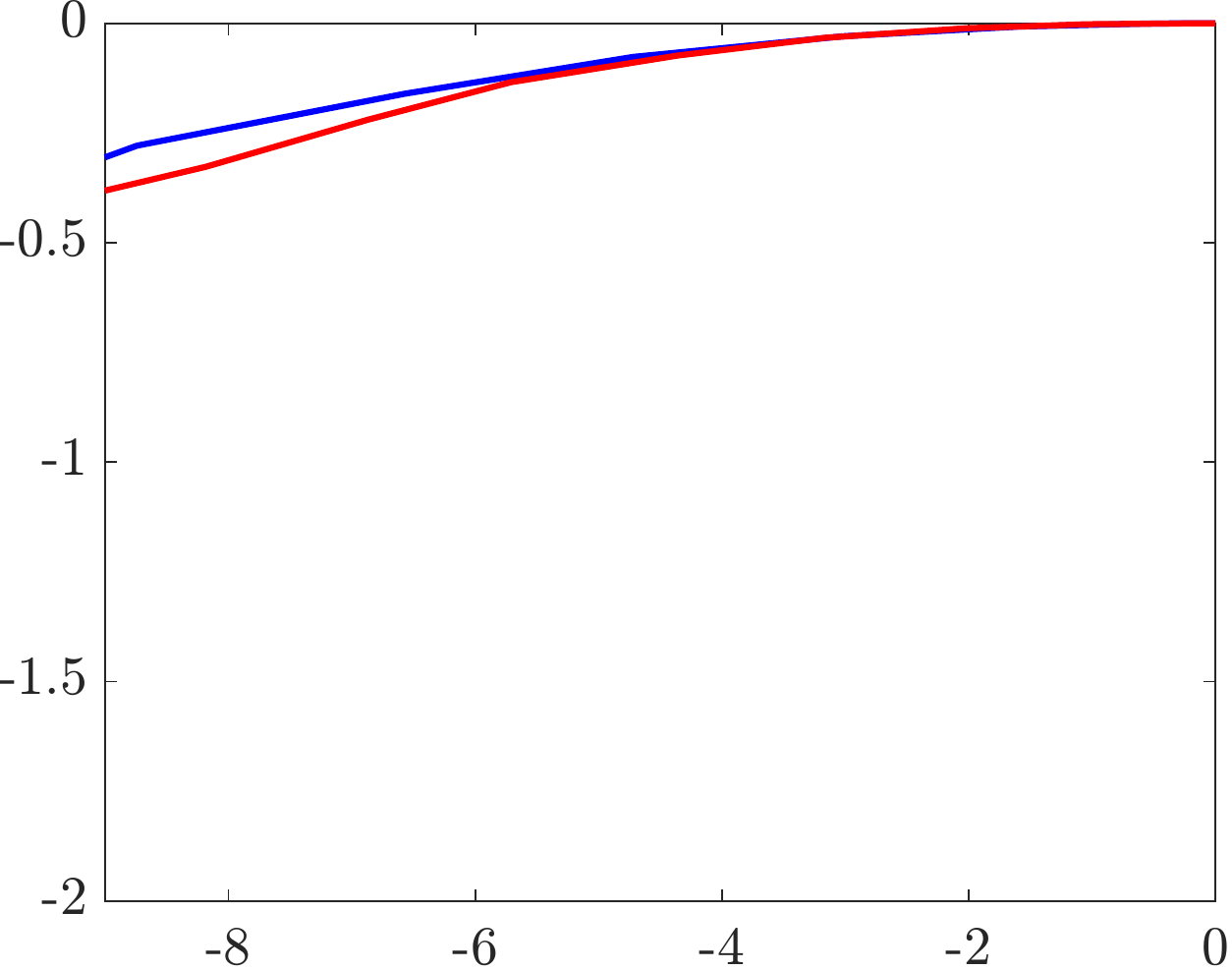}};
  \node[above of= img1, node distance=0cm, yshift=1.7cm,font=\color{black}]  {\small $N=150$};   
  \node[below of= img1, node distance=0cm, yshift=-1.7cm,font=\color{black}]  {\small $\log$ FPR};
  \node[above of= img1, node distance=0cm, xshift=-1.5cm,yshift=1.7cm,font=\color{black}]  {\small (g)};
    % \node[left of= img1, node distance=0cm, rotate=90, anchor=center,yshift=2.1cm,font=\color{black}] { \small TPR};
\end{tikzpicture}\columnbreak
\begin{tikzpicture}
  \node (img1)  {\includegraphics[width=0.20\textwidth]{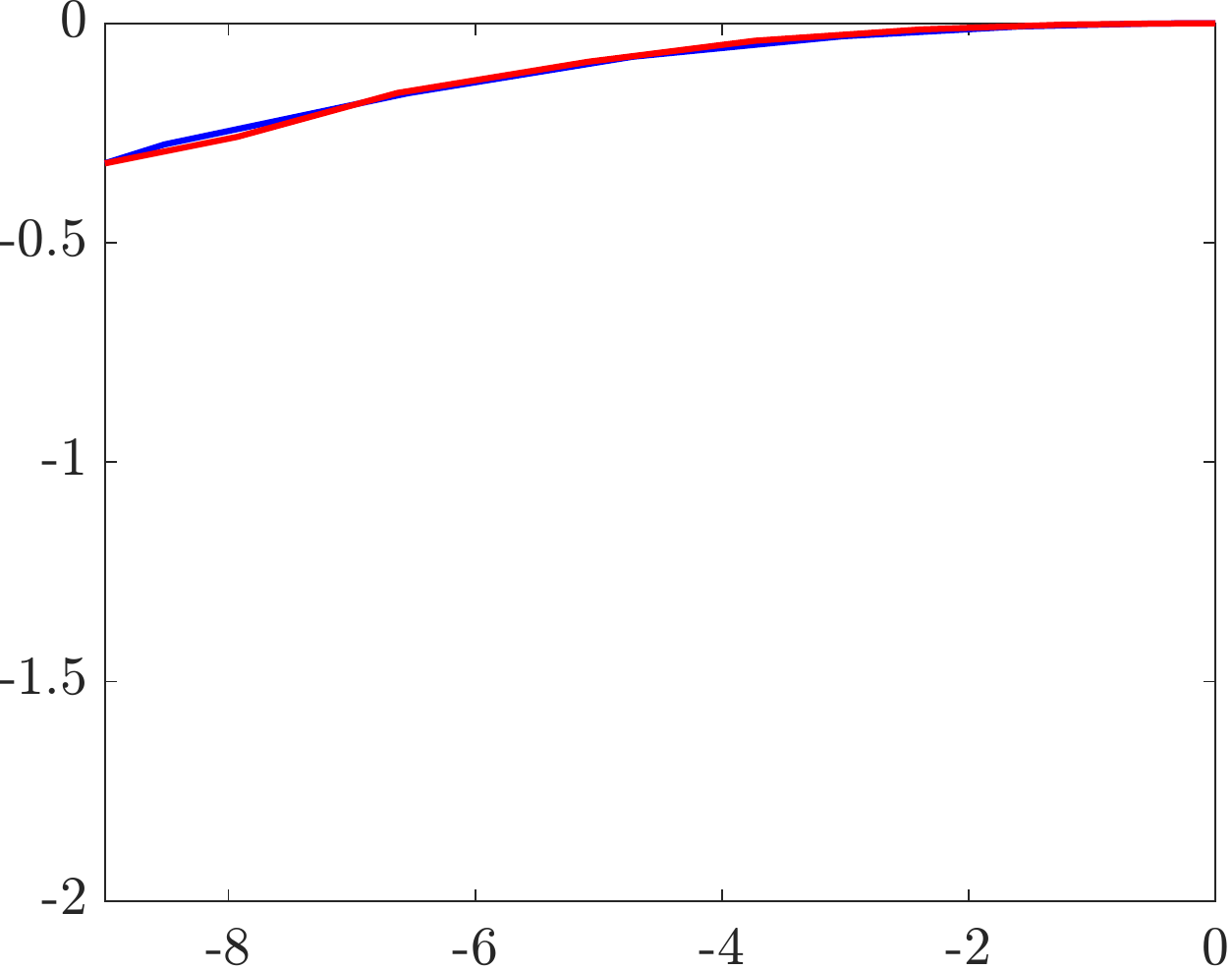}};
  \node[above of= img1, node distance=0cm, yshift=1.7cm,font=\color{black}]  {\small $N=200$};   
    \node[below of= img1, node distance=0cm, yshift=-1.7cm,font=\color{black}]  {\small $\log$ FPR};
  \node[above of= img1, node distance=0cm, xshift=-1.5cm,yshift=1.7cm,font=\color{black}]  {\small (h)};
  % \node[left of= img1, node distance=0cm, rotate=90, anchor=center,yshift=2.1cm,font=\color{black}] { \small TPR};
\end{tikzpicture}
\end{multicols}
\caption{ROC curves for the two proposed methods. For (a)-(d), detection horizon $T = 7$. The direct method performs worse than the indirect method in regimes of low data and soon outperforms the indirect method as the size of the data set $N$ increases. For (e)-(h) detection horizon $T=14$. Both methods show better performance, however the gap between the two methods increases when the data size is low and direct method takes more data to outperform the indirect method.}
\label{fig:ROCcompare-snr1}
\end{figure*}
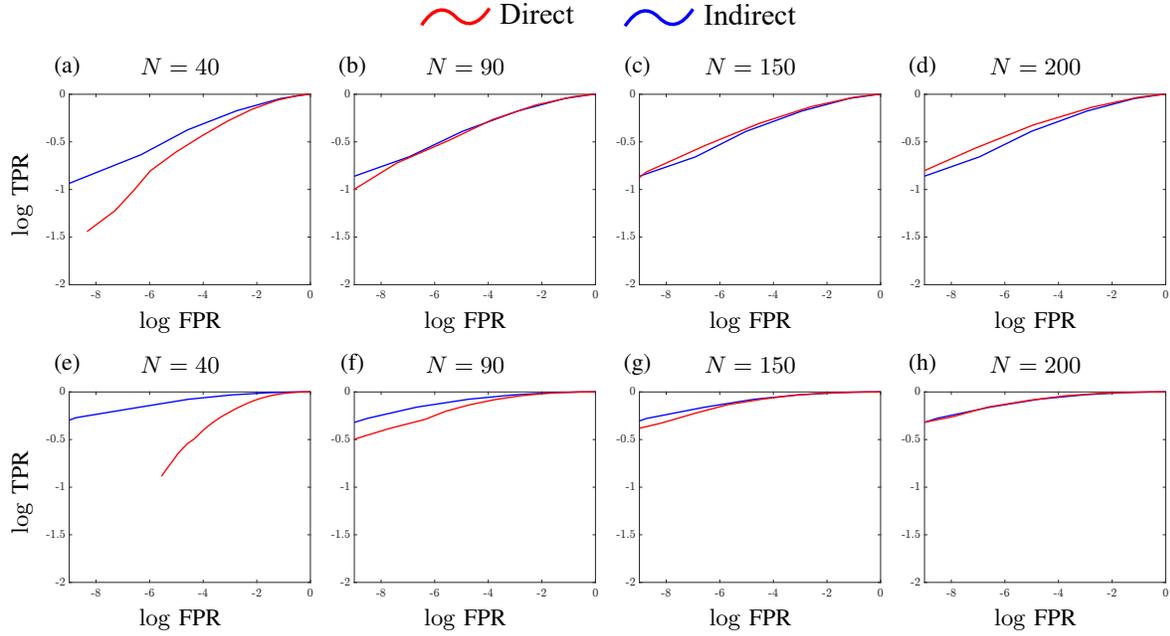

In this comparison, we set $\sigma_u = \sigma_v
 = \sigma_w = 1$ and compare the ROC curves of the two methods by plotting the
 logarithm of the FPR versus the logarithm of the TPR. In Figure \ref
 {fig:ROCcompare-snr1}, it is evident that the direct method performs poorly at
 low data-regimes while the indirect method performs better. However, as the
 data size $N$ increases, the direct method outperforms the indirect method. As
 the detection horizon increases to $T=14$ in Figure \ref{fig:ROCcompare-snr1}
 (e)-(h), the direct method takes longer to outperform the indirect method.
 However, both methods feature higher TPR as $T$ increases.

\subsubsection{Comparison 2}
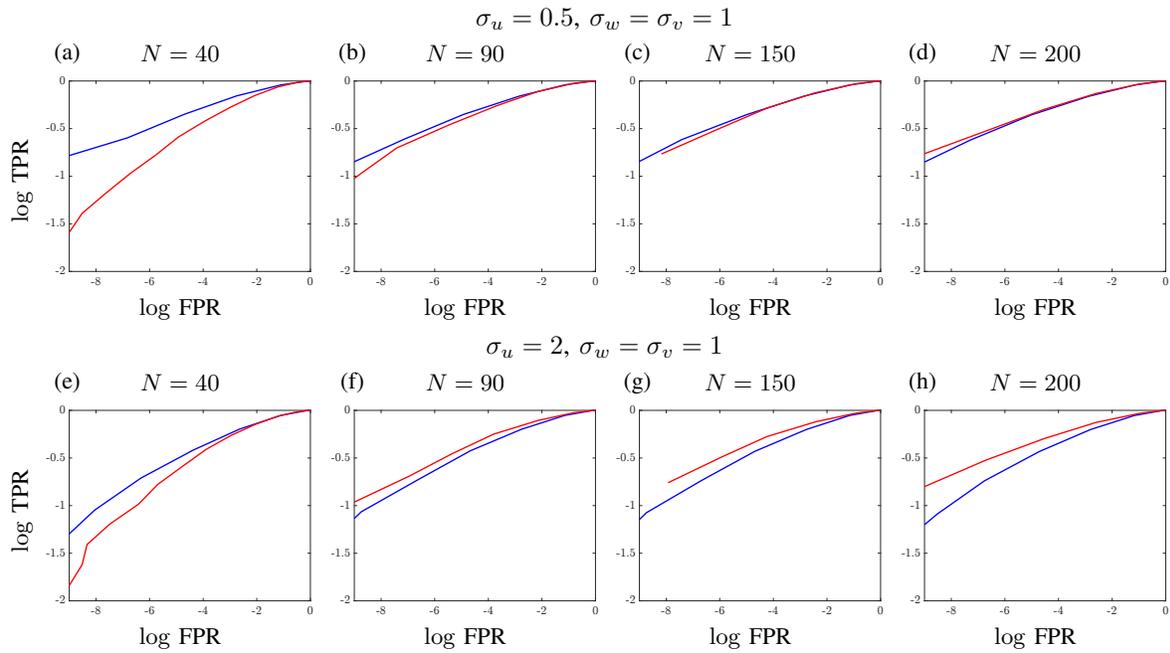
\begin{figure*}
\centering
$\sigma_u=0.5$, $\sigma_w = \sigma_v = 1$
\vspace*{-0.3cm}
\begin{multicols}{4}
\hspace*{0.75cm}
\begin{tikzpicture}
  \node (img1)  {\includegraphics[width=0.20\textwidth]{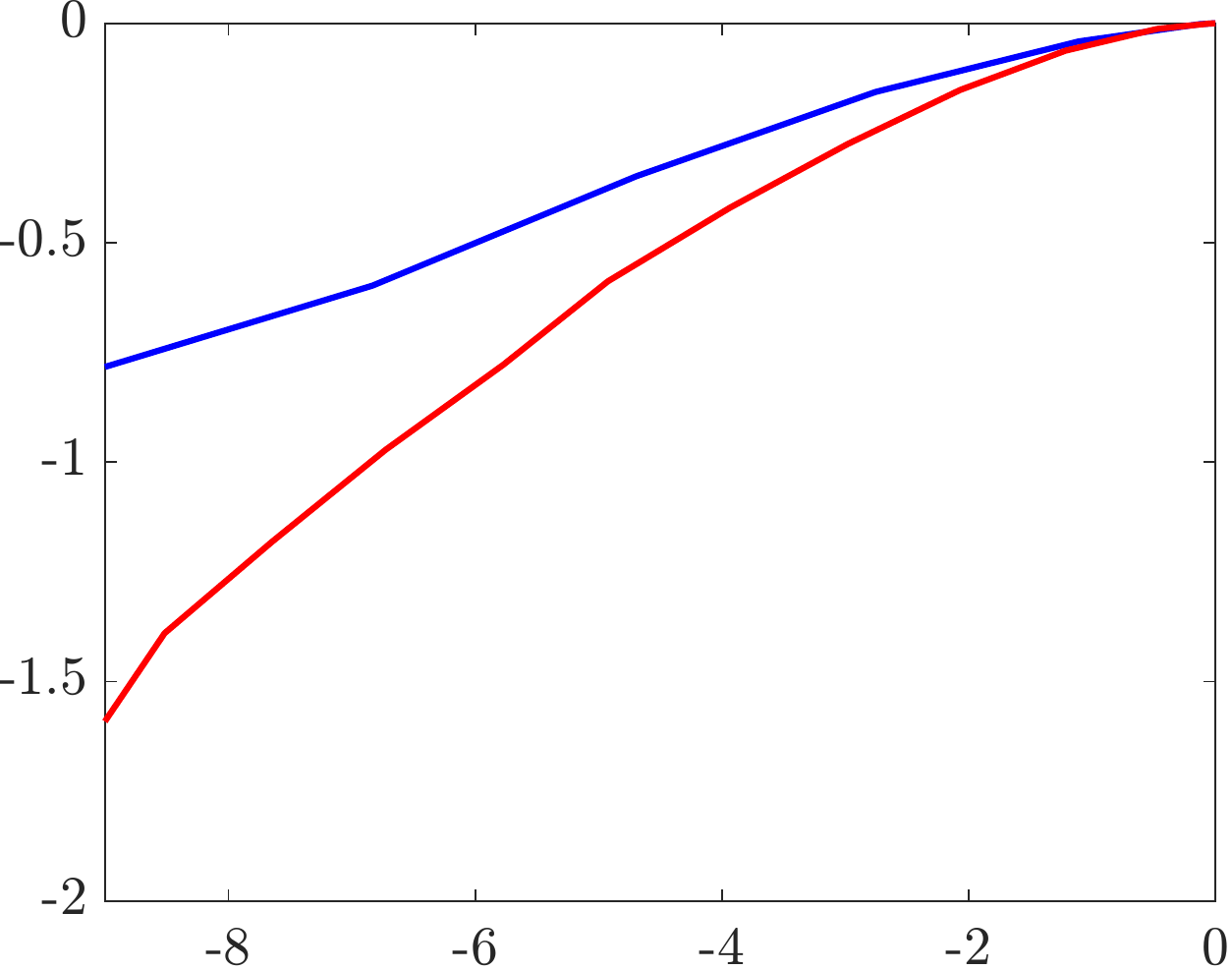}};
  \node[above of= img1, node distance=0cm, yshift=1.7cm,font=\color{black}]  {\small $N=40$};  
  \node[below of= img1, node distance=0cm, yshift=-1.7cm,font=\color{black}]  {\small $\log$ FPR};
  \node[above of= img1, node distance=0cm, xshift=-1.5cm,yshift=1.7cm,font=\color{black}]  {\small (a)};
  \node[left of= img1, node distance=0cm, rotate=90, anchor=center,yshift=2.1cm,font=\color{black}] { \small $\log$ TPR};
\end{tikzpicture}\columnbreak
\begin{tikzpicture}
  \node (img1)  {\includegraphics[width=0.20\textwidth]{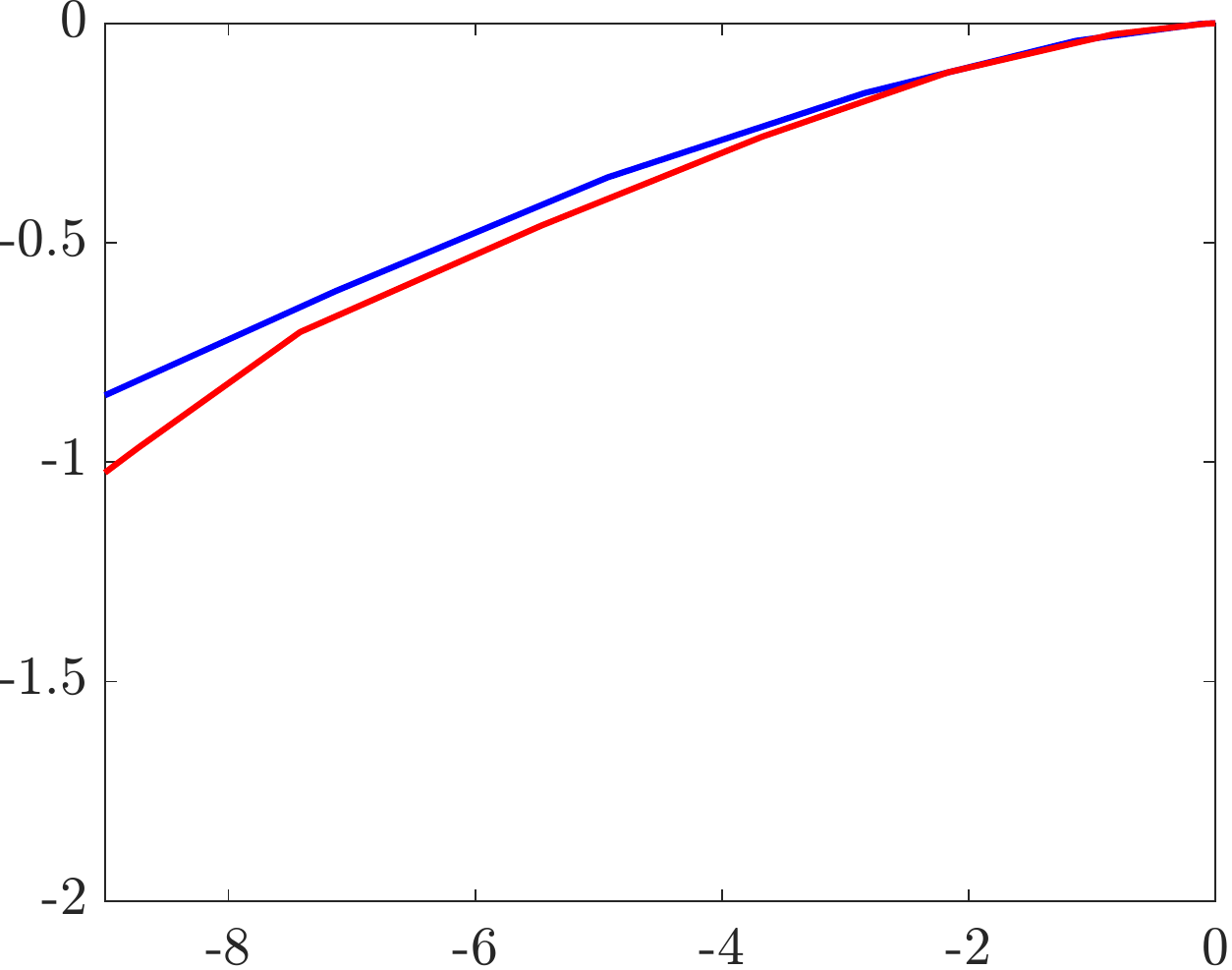}};
  \node[above of= img1, node distance=0cm, yshift=1.7cm,font=\color{black}]  {\small $N=90$}; 
  \node[below of= img1, node distance=0cm, yshift=-1.7cm,font=\color{black}]  {\small $\log$ FPR};
  \node[above of= img1, node distance=0cm, xshift=-1.5cm,yshift=1.7cm,font=\color{black}]  {\small (b)};
  % \node[left of= img1, node distance=0cm, rotate=90, anchor=center,yshift=2.1cm,font=\color{black}] { \small TPR};
\end{tikzpicture}\columnbreak 
\begin{tikzpicture}
  \node (img1)  {\includegraphics[width=0.20\textwidth]{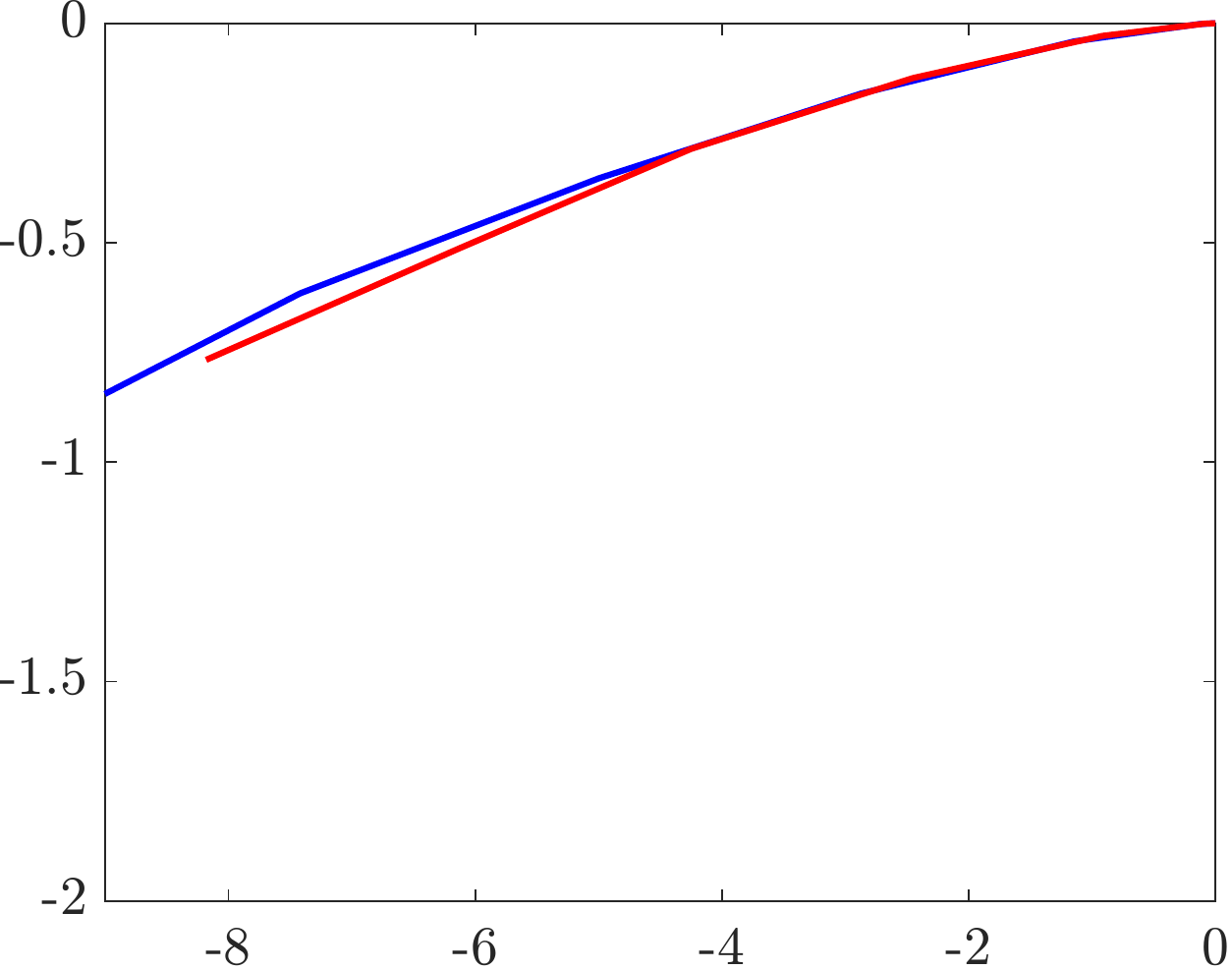}};
  \node[above of= img1, node distance=0cm, yshift=1.7cm,font=\color{black}]  {\small $N=150$};   
   \node[below of= img1, node distance=0cm, yshift=-1.7cm,font=\color{black}]  {\small $\log$ FPR};
  \node[above of= img1, node distance=0cm, xshift=-1.5cm,yshift=1.7cm,font=\color{black}]  {\small (c)};
  % \node[left of= img1, node distance=0cm, rotate=90, anchor=center,yshift=2.1cm,font=\color{black}] { \small TPR};
\end{tikzpicture}\columnbreak
\begin{tikzpicture}
  \node (img1)  {\includegraphics[width=0.20\textwidth]{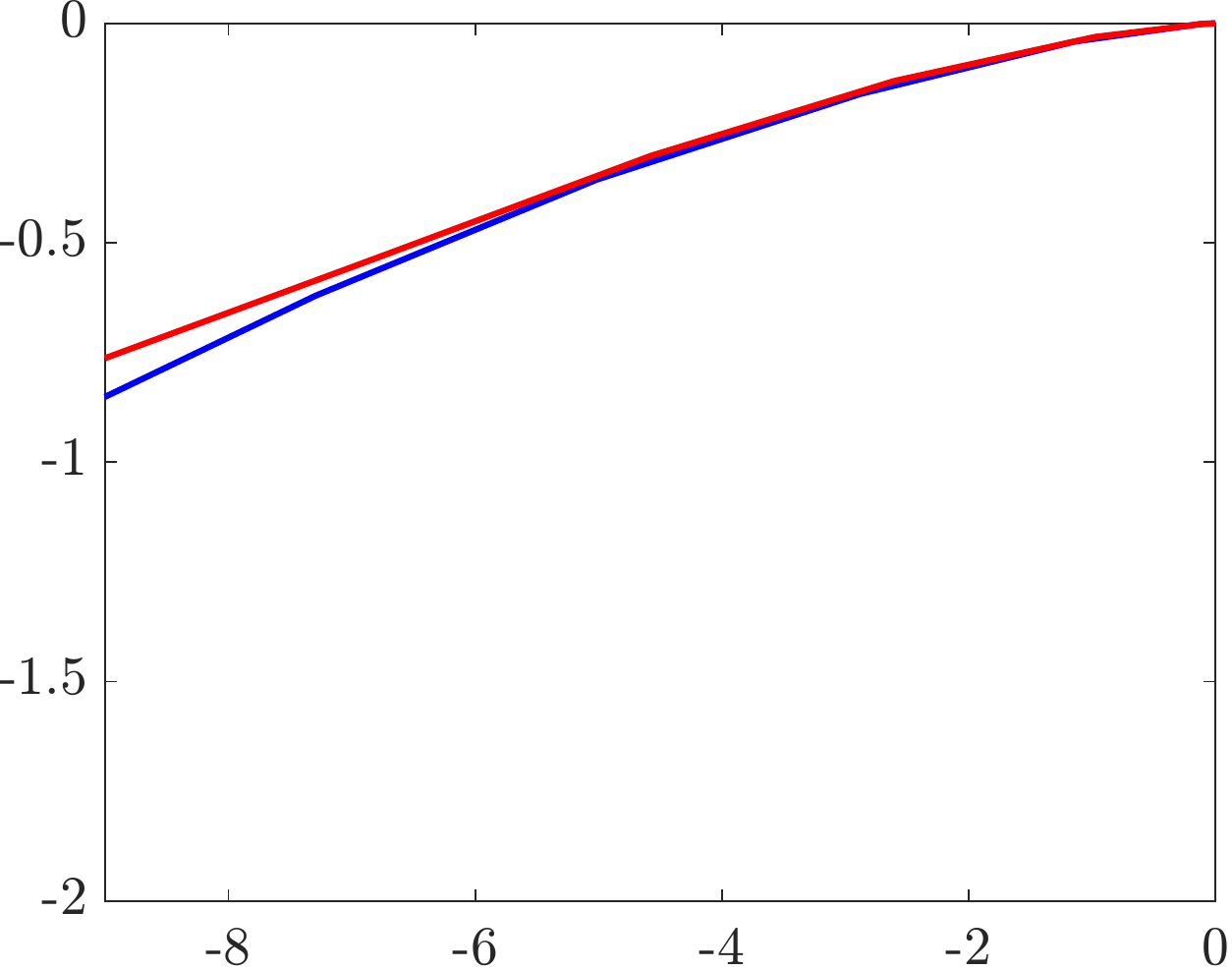}};
  \node[above of= img1, node distance=0cm, yshift=1.7cm,font=\color{black}]  {\small $N=200$};   
  \node[below of= img1, node distance=0cm, yshift=-1.7cm,font=\color{black}]  {\small $\log$ FPR};
  \node[above of= img1, node distance=0cm, xshift=-1.5cm,yshift=1.7cm,font=\color{black}]  {\small (d)};
  % \node[left of= img1, node distance=0cm, rotate=90, anchor=center,yshift=2.1cm,font=\color{black}] { \small TPR};
\end{tikzpicture}
\end{multicols}
\vspace*{-0.5cm}
$\sigma_u=2$, $\sigma_w = \sigma_v = 1$
\vspace*{-0.3cm}
\begin{multicols}{4}
\hspace*{0.75cm}
\begin{tikzpicture}
  \node (img1)  {\includegraphics[width=0.20\textwidth]{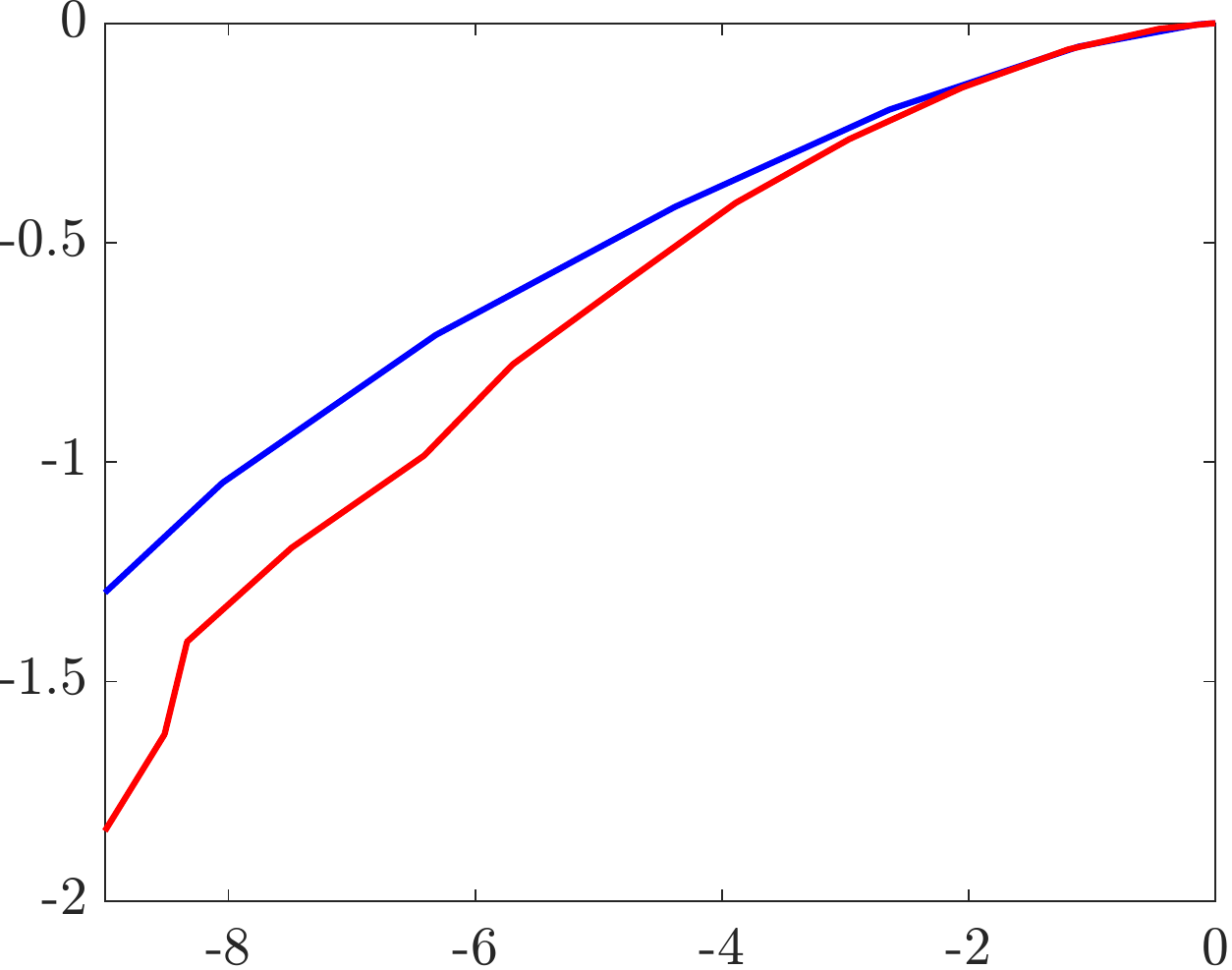}};
  \node[above of= img1, node distance=0cm, yshift=1.7cm,font=\color{black}]  {\small $N=40$};  
  \node[below of= img1, node distance=0cm, yshift=-1.7cm,font=\color{black}]  {\small $\log$ FPR};
  \node[above of= img1, node distance=0cm, xshift=-1.5cm,yshift=1.7cm,font=\color{black}]  {\small (e)};
  \node[left of= img1, node distance=0cm, rotate=90, anchor=center,yshift=2.1cm,font=\color{black}] { \small $\log$ TPR};
\end{tikzpicture}\columnbreak
\begin{tikzpicture}
  \node (img1)  {\includegraphics[width=0.20\textwidth]{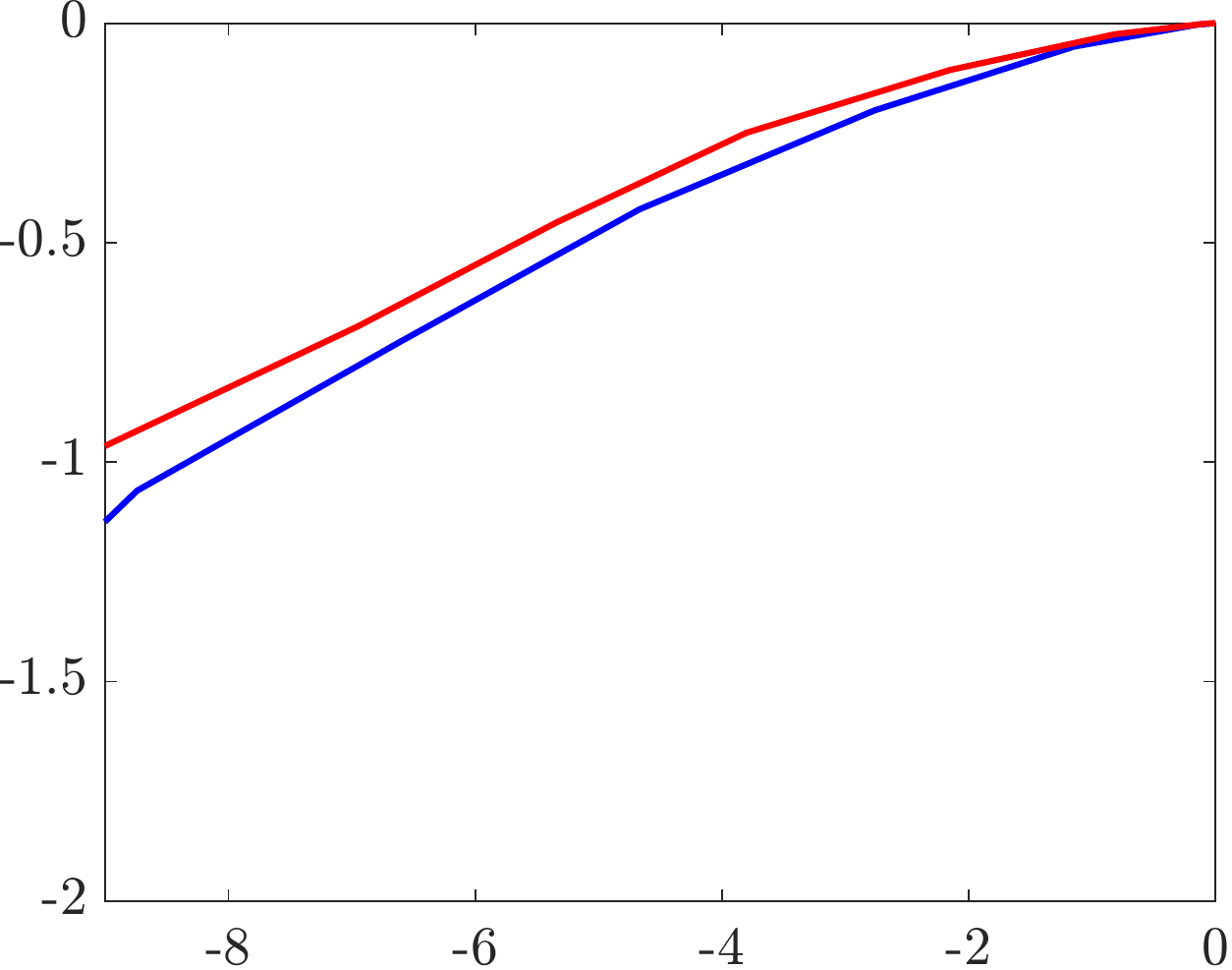}};
  \node[above of= img1, node distance=0cm, yshift=1.7cm,font=\color{black}]  {\small $N=90$}; 
  \node[below of= img1, node distance=0cm, yshift=-1.7cm,font=\color{black}]  {\small $\log$ FPR};
  \node[above of= img1, node distance=0cm, xshift=-1.5cm,yshift=1.7cm,font=\color{black}]  {\small (f)};
  % \node[left of= img1, node distance=0cm, rotate=90, anchor=center,yshift=2.1cm,font=\color{black}] { \small TPR};
\end{tikzpicture}\columnbreak 
\begin{tikzpicture}
  \node (img1)  {\includegraphics[width=0.20\textwidth]{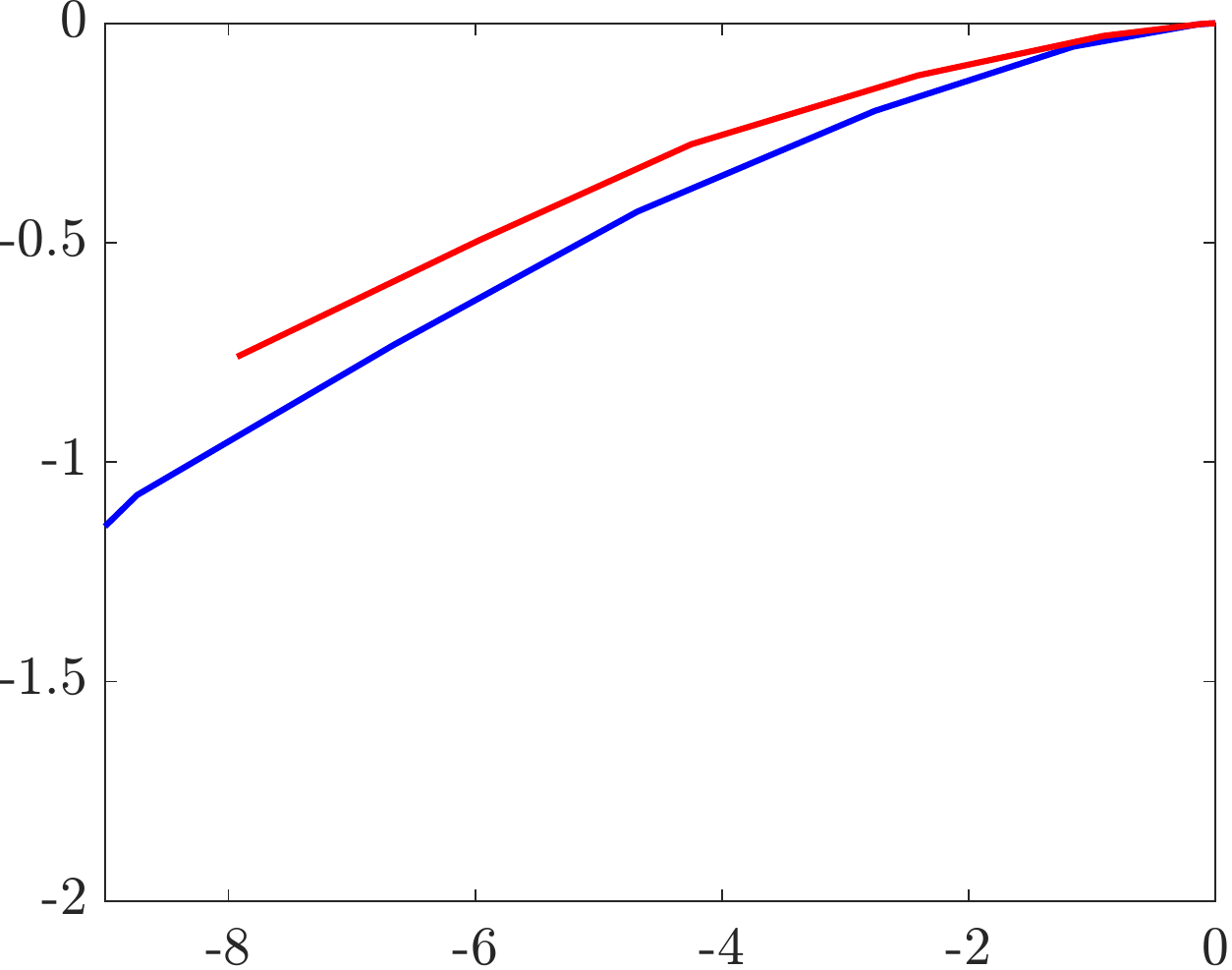}};
  \node[above of= img1, node distance=0cm, yshift=1.7cm,font=\color{black}]  {\small $N=150$};   
   \node[below of= img1, node distance=0cm, yshift=-1.7cm,font=\color{black}]  {\small  $\log$ FPR};
  \node[above of= img1, node distance=0cm, xshift=-1.5cm,yshift=1.7cm,font=\color{black}]  {\small (g)};
  % \node[left of= img1, node distance=0cm, rotate=90, anchor=center,yshift=2.1cm,font=\color{black}] { \small TPR};
\end{tikzpicture}\columnbreak
\begin{tikzpicture}
  \node (img1)  {\includegraphics[width=0.20\textwidth]{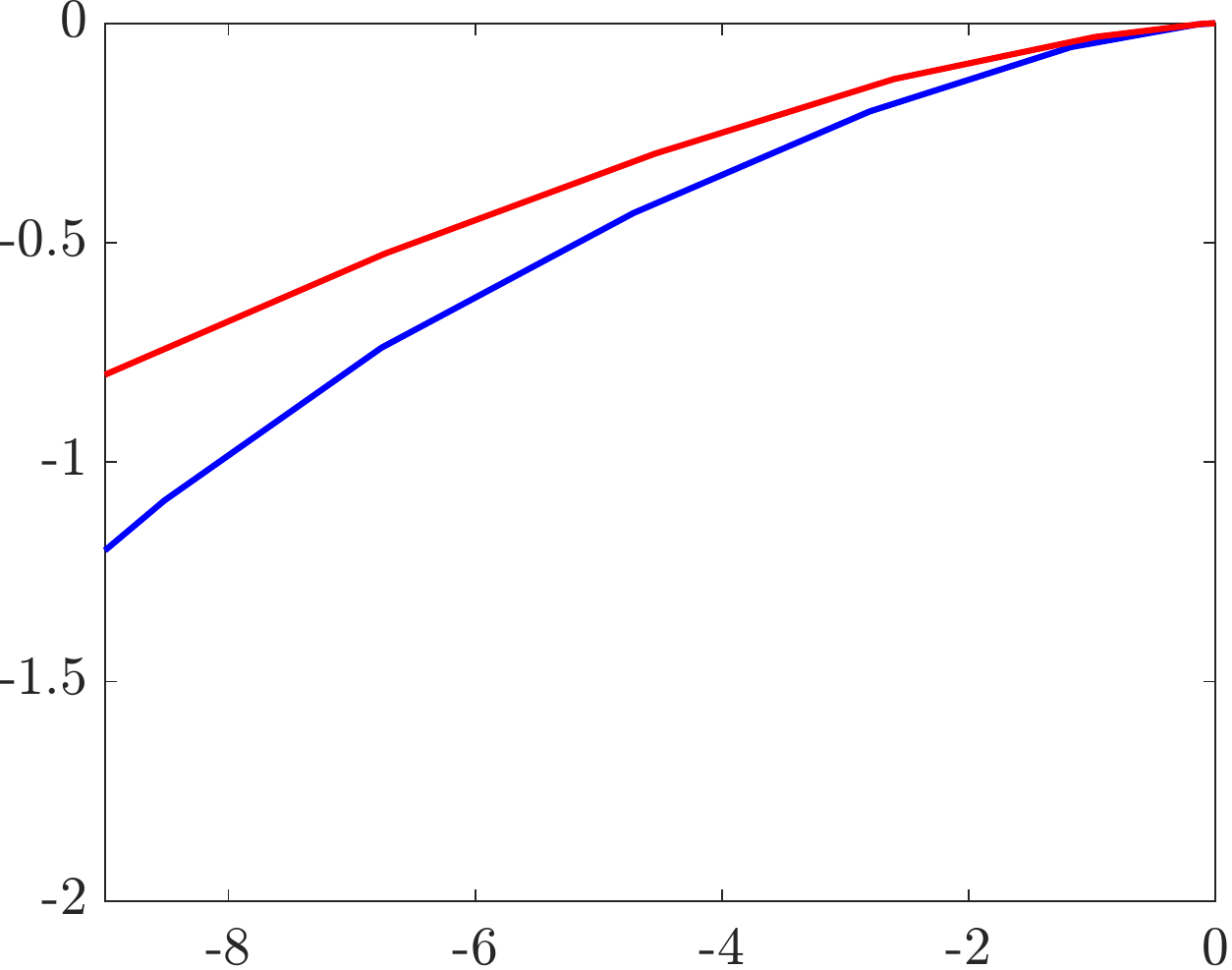}};
  \node[above of= img1, node distance=0cm, yshift=1.7cm,font=\color{black}]  {\small $N=200$};   
  \node[below of= img1, node distance=0cm, yshift=-1.7cm,font=\color{black}]  {\small $\log$ FPR};
  \node[above of= img1, node distance=0cm, xshift=-1.5cm,yshift=1.7cm,font=\color{black}]  {\small (h)};
  % \node[left of= img1, node distance=0cm, rotate=90, anchor=center,yshift=2.1cm,font=\color{black}] { \small TPR};
\end{tikzpicture}
\end{multicols}
\caption{ROC curves for different SNR with detection horizon $T=7$. For (a)-(d) $\sigma_u=0.5$, $\sigma_w = \sigma_v = 1$, which means that SNR is low. The direct method performs worse for the same amount of data with low SNR. For figures (e)-(h) $\sigma_u=2$, $\sigma_w = \sigma_v = 1$, which means that SNR is high. The direct method performs better and takes fewer samples to outperform the indirect method. However, the indirect method is affected less by the change in SNR.}
\label{fig:ROCcompare-snr0.5}
\end{figure*}

In this comparison, we vary the the signal to noise ratio in the data. Similar to the first comparison, we compare the ROC curves of the two methods. We first set $\sigma_u = 0.5$ and $\sigma_v = \sigma_w = 1$. From Figure \ref{fig:ROCcompare-snr0.5}(a)-(d), we can see that the performance of both the methods deteriorates as the SNR decreases. However, the indirect method is consistently better when the size of the data set is small. The direct method is affected more by the decrease in SNR with a decrease in TPR. Further, it takes more data for the direct method to outperform the indirect method. Next, we set $\sigma_u = 2$ and $\sigma_v = \sigma_w = 1$, thereby increasing the SNR. In Figure \ref{fig:ROCcompare-snr0.5}(e)-(h), the gap between the direct and indirect methods decreases across all data sizes. Further, with more data the direct method outperforms the indirect method quicker. 

\section{Conclusion and future work}

In this paper, we proposed a data-driven $\chi^2$ attack detector that uses the
behaviors of the system to differentiate between nominal and attacked
operation. The proposed detector requires estimation of the covariance of
input-output behaviors of the system. This is achieved by the two proposed methods,
the direct and indirect method. We analytically showed the consistency of the
two methods and established a probabilistic finite-sample error bounds. After a
numerical study of the performance of the two methods, it was evident that the
neither method is invariably superior. The direct method performs well in cases
with more data and shorter detection horizons. However, as the attack detection
horizon increases, the direct method starts to perform poorly. The indirect
method, with its reliance on an underlying generative model of the system, 
outperforms the direct method in small data regimes. This study
shows that neither method has invariably superior performance and that the choice of
method must be based on the size and characteristics of the available
data set.

\bibliographystyle{unsrt}

\bibliography{refs,alias,FP,Main,New}

\begin{thebibliography}{10}

\bibitem{SMD-MP-etal:2019}
S.~M. Dibaji, M.~Pirani, D.~B. Flamholz, A.~M. Annaswamy, K.~H. Johansson, and
  A.~Chakrabortty.
\newblock A systems and control perspective of cps security.
\newblock {\em Annual reviews in control}, 47:394--411, 2019.

\bibitem{FP-FD-FB:10y}
F.~Pasqualetti, F.~D{\"o}rfler, and F.~Bullo.
\newblock Attack detection and identification in cyber-physical systems.
\newblock {\em IEEE Transactions on Automatic Control}, 58(11):2715--2729,
  2013.

\bibitem{CZB-FP-VG:2017}
C.-Z. Bai, F.~Pasqualetti, and V.~Gupta.
\newblock Data-injection attacks in stochastic control systems: Detectability
  and performance tradeoffs.
\newblock {\em Automatica}, 82:251--260, 2017.

\bibitem{YM-RC-BS:2013}
Y.~Mo, R.~Chabukswar, and B.~Sinopoli.
\newblock Detecting integrity attacks on scada systems.
\newblock {\em IEEE Transactions on Control Systems Technology},
  22(4):1396--1407, 2013.

\bibitem{CM-JR:2016}
C.~Murguia and J.~Ruths.
\newblock Cusum and chi-squared attack detection of compromised sensors.
\newblock In {\em 2016 IEEE Conference on Control Applications (CCA)}, pages
  474--480, Buenos Aires, Argentina, September 2016.

\bibitem{RT-CM-JR:2018}
R.~Tunga, C.~Murguia, and J.~Ruths.
\newblock Tuning windowed chi-squared detectors for sensor attacks.
\newblock In {\em 2018 Annual American Control Conference (ACC)}, pages
  1752--1757, Milwaukee, USA, June 2018.

\bibitem{YL-LS-TC:2017}
Y.~Li, L.~Shi, and T.~Chen.
\newblock Detection against linear deception attacks on multi-sensor remote
  state estimation.
\newblock {\em IEEE Transactions on Control of Network Systems}, 5(3):846--856,
  2017.

\bibitem{AYL-GHY:2022}
A.-Y. Lu and G.-H. Yang.
\newblock False data injection attacks against state estimation without
  knowledge of estimators.
\newblock {\em IEEE Transactions on Automatic Control}, 2022.

\bibitem{JG-SA-MT-ZSL:17}
J.~Goh, S.~Adepu, M.~Tan, and Z.~S. Lee.
\newblock Anomaly detection in cyber physical systems using recurrent neural
  networks.
\newblock In {\em {IEEE} 18th {I}nternational {S}ymposium on {H}igh {A}ssurance
  {S}ystems {E}ngineering (HASE)}, pages 140--145, 2017.

\bibitem{MK-AS:2021}
M.~Kravchik and A.~Shabtai.
\newblock Efficient cyber attack detection in industrial control systems using
  lightweight neural networks and pca.
\newblock {\em IEEE Transactions on Dependable and Secure Computing}, 2021.

\bibitem{YL-etal:2021}
Y.~Luo, Y.~Xiao, L.~Cheng, G.~Peng, and D.~Yao.
\newblock Deep learning-based anomaly detection in cyber-physical systems:
  Progress and opportunities.
\newblock {\em ACM Computing Surveys (CSUR)}, 54(5):1--36, 2021.

\bibitem{AT-etal:2014}
A.~Tiwari, B.~Dutertre, D.~Jovanovi\'{c}, T.~de~Candia, P.~D. Lincoln,
  J.~Rushby, D.~Sadigh, and S.~Seshia.
\newblock Safety envelope for security.
\newblock In {\em 3rd International Conference on High Confidence Networked
  Systems}, page 85–94, Berlin, Germany, April 2014.

\bibitem{VK-FP:2020}
V.~Krishnan and F.~Pasqualetti.
\newblock Data-driven attack detection for linear systems.
\newblock {\em IEEE Control Systems Letters}, 5(2):671--676, 2020.

\bibitem{EN-KK:2017}
E.~Naderi and K.~Khorasani.
\newblock A data-driven approach to actuator and sensor fault detection,
  isolation and estimation in discrete-time linear systems.
\newblock {\em Automatica}, 85:165--178, 2017.

\bibitem{MT-KK-IS-NM:2021}
M.~Taheri, K.~Khorasani, I.~Shames, and N.~Meskin.
\newblock Data-driven covert-attack strategies and countermeasures for
  cyber-physical systems.
\newblock In {\em 2021 IEEE 60th Conference on Decision and Control (CDC)},
  pages 4170--4175, Austin, USA, December 2021.

\bibitem{TS-AR:2019}
T.~Sarkar and A.~Rakhlin.
\newblock Near optimal finite time identification of arbitrary linear dynamical
  systems.
\newblock In {\em 36th International Conference on Machine Learning}, volume~97
  of {\em Proceedings of Machine Learning Research}, pages 5610--5618, Long
  Beach, USA, June 2019.

\bibitem{ANB-PDM-AN:2018}
A.~N. Bishop, P.~Del~Moral, and A.~Niclas.
\newblock An introduction to wishart matrix moments.
\newblock {\em Foundations and Trends in Machine Learning}, 11(2):97--218,
  2018.

\bibitem{SZ:2012}
S.~Zhu.
\newblock A short note on the tail bound of wishart distribution.
\newblock {\em arXiv preprint arXiv:1212.5860}, 2012.

\bibitem{JAT:2012}
J.~A. Tropp.
\newblock User-friendly tail bounds for sums of random matrices.
\newblock {\em Foundations of computational mathematics}, 12(4):389--434, 2012.

\bibitem{PMG-etal:1990}
P.~M. Gahinet, A.~J. Laub, C.~S. Kenney, and G.~A. Hewer.
\newblock Sensitivity of the stable discrete-time lyapunov equation.
\newblock {\em IEEE Transactions on Automatic Control}, 35(11):1209--1217,
  1990.

\end{thebibliography}

\section*{Appendix}

\subsection{Proof of Theorem \ref{thm:consistency}} \label{proof:consistency}
We provide a sketch of the proof in two parts. We first prove the consistency of the direct method using the properties of the Wishart distribution \cite{ANB-PDM-AN:2018}. Next, we prove the consistency of the indirect method using properties of the OLS solution.

\subsubsection{Direct method}
The sample covariance matrix is an unbiased estimator of the true
covariance. Since $\bm{z}^{(i)}$ are i.i.d. gaussian experiments, $\bm{Z}\bm
{Z}^\top$ in fact follows a Wishart distribution \cite{ANB-PDM-AN:2018}. By the properties of the Wishart distribution, for every element $\hat{S}^d_{ij}$, 
\begin{align*}
\mathbb{V}\mathrm{ar}\left[\hat{S}^{d}_{ij}\right] = \frac{1}{N}\left(S^2_{ij}-S_{ii}S_{jj}\right).
\end{align*}
Therefore, as $N \rightarrow \infty$, $\mathbb{V}\mathrm{ar}\left[\hat{S}^{d}_
{ij}\right] \rightarrow 0$. By applying the Chebyshev inequality for every element $\hat{S}^d_{ij}$, it follows that
the sample covariance is a consistent estimator of the covariance $S$. 
\subsubsection{Indirect method}
First, consider the
 augmented behavior $h_t = [f_t^\top | \zeta_t^\top]^\top$, where $\zeta_t$ defined as,
\begin{align*}
\zeta_t = [w_{t-L}^\top \ w_{t-L+1}^\top \ \dots \ w_{t}^\top | v_{t-L}^\top \ v_{t-L+1}^\top \ \dots \ v_{t}^\top ]^\top. 
\end{align*} 
When $L=n$, $h_t$ follows a stationary Vector Auto-Regressive(1) process such
that $h_{t+1} = \mathcal{G}h_t + e_t$, where $e_t$ is a vector of external noises. We divide the matrix $\mathcal{G}$ into block matrices $G_{11},G_{12},G_{21},$ and $G_{22}$. $G_
{11} \in \R^{L(p+m)\times L(p+m)}$ captures the dependence of future minor behaviors $f_{t+1}$ on past minor behaviors
$f_t$ which is in the Brunovsky canonical form. $G_{12}$ captures the
dependence of $f_{t+1}$ on the past noises $\zeta_t$. $G_{21}$ is a zero matrix
as the noises $\zeta_t$ do not depend on the behaviors $f_t$. Lastly, $G_
{22}$ is the matrix that captures the dependence of future noise $\zeta_
{t+1}$ on past noises $\zeta_t$. Let the covariance of the augmented
behaviors be $\E[h_t h_t^\top]=\Hc$. Then, $\Hc = \mathrm{Cov}
[\mathcal{G}h_{t}+e_t]= \G \Hc \G^\top + \Sigma_e$, which is in
the form of a Lyapunov equation. Therefore, it is sufficient to know $\G$ and
$\Sigma_e$ to compute $\Hc$. However, we are only interested in the first $L
(m+p)\times L(m+p)$ block of $\Hc$ as it captures the covariance of the minor
behaviors $f_t$ given by $P$ from equation \eqref{eqn:indirect-lyap}. We now
describe how $P$ can be computed from the Lyapunov equation $\Hc
= \G \Hc \G^\top + \Sigma_e$.

 Define the $(1,1), (1,2), (2,1)$ and $(2,2)$ blocks of $\Hc$ as $P_{11},P_
 {12},P_{21}$ and $P_{22}$ of appropriate size. It is crucial to note that $P_{11} = P$ in equation \eqref{eqn:indirect-lyap} as it captures the covariance of the minor behaviors. Then $P_{11}$ can be computed by solving the Lyapunov equation as,
 \begin{align}
 P &= (\G_{11}+\G_{12}P_{21}P^{-1})P(\G_{11}+\G_{12}P_{21}P^{-1})^\top + \nonumber
 \\ & \ \ \ \ {\Sigma_e}_{11} + \G_{12} (\Hc / P) \G_{12}^\top. \label{eqn:indirectprooflyap}
 \end{align}
Let $\tilde
{\G} = \G_{11}+\G_{12}P_{21}P^{-1}$ and $\tilde{\Sigma} = {\Sigma_e}_
{11} + \G_{12} (\Hc / P) \G_{12}^\top$. Then equation \eqref
{eqn:indirectprooflyap} can be rewritten as $P = \tilde{\G}P\tilde
{\G} + \tilde{\Sigma}$. This implies that there exists a VAR(1) system such
that $f_{t+1} = \tilde{\G}f_t + h_{t}$, where $h_t \sim \mathcal{N}(0,\tilde
{\Sigma})$.

Next, we prove that regressing future minor behaviors over past minor behaviors
using OLS yields $\tilde{G}$. We seek to solve the following problem.
\begin{align}
\min_{\hat{G}} \  \E\left[\|f_{t+1} - \hat{G}f_t\|^2\right] \label{eqn:indirectOLSproblem}
\end{align}
Note that $f_{t+1} = \G_{11}f_t+G_{12}\zeta + \tilde{e}_t$, where $\tilde{e}_t$ is the first $L(m+p)$ elements of the vector $e_t$. Then,
\begin{align*}
&\E[\|f_{t+1} - \hat{G}f_t\|] = \E \left[\|\left((\G_{11} - \hat{G})f_t+\G_{12}\zeta_t+\tilde{e}_t\right)\|\right] \\
&=\E\Big[\left((\G_{11} - \hat{G})f_t+\G_{12}\zeta_t+\tilde{e}_t\right)^\top \\ & \quad \quad \quad \quad \quad \quad \left( (\G_{11} - \hat{G})f_t+\G_{12}\zeta_t+\tilde{e}_t\right)\Big] \\
&= \E\Big[\mathrm{Tr}\Big(\left((\G_{11} - \hat{G})f_t+\G_{12}\zeta_t+\tilde{e}_t\right)\\
 & \quad \quad \quad \quad \quad \quad \left((\G_{11} - \hat{G})f_t+\G_{12}\zeta_t+\tilde{e}_t\right)^\top \Big)\Big] \\
&= \mathrm{Tr} \Big( (\G_{11} - \hat{G}) P (\G_{11} - \hat{G})^\top \\ & \quad \quad \quad \quad \quad \quad + (\G_{11}-\hat{G}) P_{12} \G_{12}^\top + \G_{12} P_{21} (\G_{11}-\hat{G})^\top \Big)
\end{align*}
Optimality condition for the problem occurs when the derivative of the cost
function with respect to $\hat{G}$ is $0$. Differentiating the cost function
using above expression with respect to $\hat{G}$ and setting it to $0$, we
obtain,
\begin{align*}
\frac{d \E\left[\|f_{t+1} - \hat{G}f_t\|^2\right] }{d\hat{G}} &= (\G_{11} - \hat{G})P - \G_{12} P_{21} = 0 \\
\hat{G} &= \G_{12}P_{21}P^{-1}+\G_{11} = \tilde{G}.
\end{align*}
Therefore, as $N \rightarrow \infty$,
$\hat{G} \rightarrow \tilde{G}$. Residuals of the least squares regression are
now $f_{t+1} - \tilde{G}f_t$. By computing the covariance of the residuals, we
obtain 
\begin{align*}
\E\left[(f_{t+1} - \tilde{G}f_t)(f_{t+1} - \tilde{G}f_t)^\top\right] = {\Sigma_e}_{11} + \G_{12}(\Hc/P)\G_{12}^\top.
\end{align*}

Therefore the ordinary least squares problem is consistent with equation \eqref
{eqn:indirectprooflyap}. Further, from equation \eqref{eqn:indirect-lyap} we have that, $\Sigma_\epsilon = {\Sigma_e}_{11} + \G_{12}(\Hc/P)G_{12}^\top = \tilde{\Sigma}$. This also proves that the residuals are uncorrelated
with the regressors $f_t$. \hfill $\blacksquare$

\subsection{Proof of Theorem \ref{thm:finitesample-direct}}

We provide a sketch of the proof. We follow the
approach from \cite{SZ:2012}. Firstly, the i.i.d. experiments $\bm{z}^{(i)}$ satisy the Bernstein moment condition. The Bernstein moment condition states that, for any sequence of random variables $ \xi_i \sim \mathcal{N}(0,\Sigma_\xi)$, and for any $l \geq 2$, and $H>0$
\begin{align*}
\E\left[(\xi_i\xi_i^\top)^l\right] \preceq \frac{l!}{2}H^{l-2}J,
\end{align*}
where $J$ is a positive definite matrix. Then it follows from Theorem 3.6 in \cite{JAT:2012} that, for any $\theta > 0$,
\begin{align}
\lambda_1(\sum_i \xi_i\xi_i^\top) &\geq N \lambda_1(S) + \sqrt{2N\theta\lambda_1(J)} + \theta H, \label{eqn:bern1}\\
\lambda_{d}(\sum_i \xi_i\xi_i^\top) &\leq N \lambda_{d}(S) - \sqrt{2N\theta\lambda_{d}(J)} + \theta H, \label{eqn:bern2}
\end{align}
with probability not more than $d e^{-\theta}$. Here $\lambda_1$ and $\lambda_d$ denote the largest and smallest eigenvalues, respectively, $d$ is the length of the vector $\xi_i$ and $N$ is the sample size. Next, to provide error bounds we follow the Bernstein moment condition for $(\bm{z}^{(i)}{\bm{z}^{(i)}}^\top-S)$. Using Lemma 4 of \cite{SZ:2012} we get: 
\begin{align*}
\E\left[\bm{z}^{(i)}{\bm{z}^{(i)}}^\top-S)^l\right] \preceq \frac{l!}{2}H^{l-2}J, \\
\E\left[S-\bm{z}^{(i)}{\bm{z}^{(i)}}^\top\right] \preceq \frac{l!}{2}H^{l-2}J,
\end{align*}
where $J = \mathrm{Tr}(S)S$ and $H=2\mathrm{Tr}(S)$. Using $r = \frac{Tr(S)}{S}$, we have that $\|J\|\leq(r+1)\|S\|^2$ and $H = 2r\|S\|$. From equations \eqref{eqn:bern1} and \eqref{eqn:bern2}, we get for any $\theta > 0$,
\begin{align}
\|\frac{1}{N}\bm{Z}\bm{Z}^\top - S \| \geq \left(\sqrt{\frac{2\theta(r+1)}{N}} + \frac{2\theta r}{N}\right)\|S\| \label{eqn:app-lowerb}
\end{align} 
with probability not greater than $2T(p+m)e^{-\theta}$. Note that equation \eqref{eqn:app-lowerb} is a lower bound on $\|\frac{1}{N}\bm{Z}\bm{Z}^\top - S \|$. The upper bound on $\|\frac{1}{N}\bm{Z}\bm{Z}^\top - S \|$ follows with probability $1-2T(p+m)e^{-\theta}$.   \hfill $\blacksquare$

\subsection{Proof of Lemma \ref{lem:fp-bnd}} 
We provide a sketch of the proof for brevity.

\emph{(a) Sensitivity of $P$:} We obtain this result from the characterizing the sensitivity of the
 Lyapunov equation \eqref{eqn:indirect-lyap}. We use the result from Corollary
 2.7 of \cite{PMG-etal:1990}

\emph{(b) Sensitivity of $\F$:}
We exploit the structure of $\hat{\F}$ to obtain the bound. 
We have the same elements in the off-diagonals of the matrix.
\begin{align*}
&\hat{\F}= 
\begin{bmatrix}
I & \hat{\M} & \dots & \hat{\M}^{T-L} \\
\hat{\M} & I & \dots & \hat{\M}^{T-L-1} \\
\vdots &\vdots &\ddots &\vdots \\
\hat{\M}^{T-L} & \hat{\M}^{T-L-1} & \dots & I
\end{bmatrix} = I\\
& + \hat{\M} \otimes \left(\begin{bmatrix}
0 & I & 0 &\dots & 0 \\
0 & 0 & I &\dots & 0 \\
\vdots &\vdots & \vdots &\ddots &\vdots \\
0 & 0 & 0 &\dots & I \\
0 & 0 & 0 &\dots & 0
\end{bmatrix} + \begin{bmatrix}
0 & I & 0 &\dots & 0 \\
0 & 0 & I &\dots & 0 \\
\vdots &\vdots & \vdots &\ddots &\vdots \\
0 & 0 & 0 &\dots & I \\
0 & 0 & 0 &\dots & 0
\end{bmatrix}^\top\right) \\
&+ \quad \quad \quad \quad \quad \vdots \\
&+ \hat{\M}^{T-L} \otimes \left(\begin{bmatrix}
0 & \dots & 0 & I \\
0 & \dots & 0 & 0 \\
\vdots & \vdots & \ddots &\vdots \\
0 & \hdotsfor{2} & 0
\end{bmatrix} + \begin{bmatrix}
0 & \dots & 0 & I \\
0 & \dots & 0 & 0 \\
\vdots & \vdots & \ddots &\vdots \\
0 & \hdotsfor{2} & 0
\end{bmatrix}^\top\right).
\end{align*}
Here, all the matrices with only the identity matrices on the off diagonal positions have an operator norm equal to 1 as they have utmost 1 identity per column. We use the property that for any arbitrary matrices $X$ and $Y$ of appropriate dimensions, $\|X \otimes Y\| = \|X\| \ \| Y\|$. Consequently,
\begin{align*}
\|\hat{\F}\| & \leq \| 1 + 2 \|\hat{\M}\| + 2 \|\hat{\M}^2\| + \dots + 2 \|\hat{\M}^{T-L}\| \ \| \\
&\leq 1 + 2 \|\hat{\M}\| + 2 \|\hat{\M}\|^2 + \dots + 2 \|\hat{\M}\|^{T-L} \quad \quad \quad \ \blacksquare
\end{align*}
\end{document}